\documentclass[10pt, indent]{article}

\usepackage{amsmath}
\usepackage{amsthm}
\usepackage{amssymb}
\usepackage{indentfirst}
\usepackage{booktabs}
\usepackage{graphicx}

\usepackage {tikz}
\usetikzlibrary {positioning}
\usepackage[obeyspaces]{url}
\usepackage{hyperref}

\newtheorem{theorem}{Theorem}

\theoremstyle{definition}

\newcommand{\al}{\alpha}
\newcommand{\be}{\beta}
\newcommand{\ga}{\gamma}
\newcommand{\ep}{\epsilon}
\newcommand{\la}{\lambda}
\newcommand{\si}{\sigma}

\newcommand{\om}{\omega}

\newcommand{\nin}{\noindent}
\newcommand{\tb}{\textbf}
\newcommand{\ti}{\textit}

\usepackage[justification=centering]{caption}

\usepackage[OT2,T1]{fontenc}
\usepackage[russian,english]{babel}
\usepackage{hyperref}

\begin{document}
\title{In praise (and search) of J.\ V.\ Uspensky}
\author{ Persi Diaconis\thanks{Research of Persi Diaconis supported by the NSF under grant number 1954042.} \\
Stanford University \and
Sandy Zabell\\
Northwestern University}
\date{}
\maketitle

\begin{abstract}
The two of us have shared a fascination with James Victor Uspensky's 1937 textbook \ti{Introduction to Mathematical Probability} ever since our graduate student days: it contains many interesting results not found in other books on the same subject in the English language, together with many non-trivial examples, all clearly stated with careful proofs. We present some of Uspensky's gems to a modern audience hoping to tempt others to read Uspensky for themselves, as well as report on a few of the other mathematical topics he also wrote about (for example, his book on number theory contains early results about perfect shuffles).

Uspensky led an interesting life: a member of the Russian Academy of Sciences, he spoke at the 1924 International Congress of Mathematicians in Toronto before leaving Russia in 1929 and coming to the US and Stanford. Comparatively little has been written about him in English; the second half of this paper attempts to remedy this.
\end{abstract}

\section{Introduction}

In 1927, when Harald Cram\'{e}r visited England and mentioned to G. H. Hardy (his friend and former teacher) he had become interested in probability theory, Hardy replied ``there was no mathematically satisfactory book in English on this subject, and encouraged me to write one'' (Cram\'{e}r, 1976, p.\ 516).  Ten years later J.\ V.\ Uspensky's book \ti{Introduction to Mathematical Probability} filled this vacuum:  ``his mathematically demanding book was the standard text on probability theory in [the US] until the appearance'' of Feller's 1950 classic, \ti{Introduction to Probability Theory and its Applications} (Reeds et al., 2015, p.\ 588).\footnote{See also Seneta, 2006, p.\ 6 (``Uspensky's book seems to have brought analytical probability, in the St. Petersburg tradition, to the United States, where it remained a primary probabilistic source until the appearance of W.\ Feller's \ti{An Introduction to Probability Theory and Its Applications} in 1951''); Bhatia, 2008, p.\ 26, quoting S.\ R.\ S.\ Varadhan (``There were not too many books available at that time [the mid-1950s]. Feller's book had just come out. Before that there was a book by Uspensky. These were the only books on Probability'').}$^{,}$\footnote{Not surprisingly, Uspensky also became a major source for less mathematically demanding textbooks such as M.\ E.\ Munroe's  \ti{Theory of Probability}, which said of Uspensky ``we recommend again and again as collateral reading'' (Munroe, 1951, p. 2).}

Who was the author of this important book?  J.\ V.\ (James Victor) Uspensky 
(born Yakov Viktorovich Uspensky, {\cyr Yakov Viktorovich Uspenski{\u i}}, 1883--1947) entered the Saint Petersburg Imperial University in 1903, received his undergraduate degree there in 1906 as well as his graduate degree in 1911, studying under the great Andrei Andreievich Markov.  In 1915 he became a Professor at Petrograd (the new name for Saint Petersburg) University, and in 1921 elected a member of the Russian Academy of Sciences.  He gave a talk at the 1924 International Congress of Mathematicians, and traveled to the United States several times during this period:  1924, 1926--1927 (when he taught at Carleton College in Minnesota), and 1929, when he moved permanently to the US and Stanford.  Initially appointed an acting Professor of Mathematics (1929--1931), he became a permanent member of Stanford's mathematics department in 1931 and remained there until his death in 1947.  He is best known today for three textbooks he wrote during his time at Stanford although, as we shall see, he was much more active than just this.  

In the following we focus on Uspensky's three textbooks, primarily the one on mathematical probability, but then circle back and take a closer look at this very interesting person.  While this paper is intended primarily as a contribution to the history of our subject, it is worth mentioning that we have used some of the  topics discussed below as  course material in both introductory and graduate courses. Historical material often brings a course to life for students. The current proofs in textbooks are frequently so streamlined it can be easy to forget just how challenging the initial proofs were. Some of the material in the sections below can be used directly as part of a lecture, or  the basis of a challenging problem set, or as reading in conjunction with a student project.\footnote{We thank Alexei Borodin and Ilya Khayutin for help with the Russian literature, Stanford Librarian Ashley Jester for help with the Stanford Archives,  Sunny Scott and her Stanford history experts for information about Uspensky during his time at Stanford, and  Stephen DeSalvo, Stewart Ethier, Jim Pitman, Eugenio Regazzini, Reinhard Siegmund-Schultze, and Steve Stigler for helpful comments.}      

\section{Introduction to Mathematical Probability}

One of the great strengths of Uspensky's book was its making available in English for the first time a substantial fraction of the Russian literature on probability.  In the following we discuss a few particularly interesting examples.

\subsection{Computing binomial tail probabilities}

Suppose $S_n$ is the number of heads in $n$ independent tosses of a $p$-coin.  If $n = 9,000$ and $p = 1/3$,  the normal approximation tells us that $P(S_{9000} > 3090)$ is about 0.02209.  R tells us the exact value to five places is 0.02170, so the normal approximation here is accurate to three places.  But how can we tell this without R?  Estimating the error in the normal approximation or doing better by adding correction terms to it can be tricky.  It turns out there is a very clever way of tackling the twin issues of better approximation and error estimation without appealing at all to the central limit theorem (or brute force calculation).  

In Chapter 3, Section 8 of his book, Uspensky notes that although there are asymptotic formulas for estimating interval binomial probabilities of the form $P(a < S_n < b)$, ``less known is the ingenious use by Markoff of continued fractions for that purpose''.    In the following we give (1) Markov's method as described by Uspensky, together with a proof of the key convergent inequalities;   (2) a discussion of the method, including an important forward recursion the convergents satisfy permitting their efficient computation; and (3) comment briefly on some later literature.  

\subsubsection{Markov's method of continued fractions}

Let $l$ an integer greater than $np$, and consider the right binomial tail probability
\[
P(l) := P(S_n > l)
\]
Then a trite calculation yields
\[
P(l) = b(l + 1; n , p) 
\left[ 1 + \frac{n - l -1 }{l + 2} \frac{p}{q} + \frac{(n - l - 1)(n - l- 2)}{(l + 2)(l + 3)} \left(\frac p q \right)^2 + \dots \right].
\]
The term in square brackets is a special case of the hypergeometric series
\[
F(\al, \be, \ga, x) = 1 + \frac{\al \be}{1 \cdot \ga} x + \frac{\al (\al + 1)\be(\be + 1)}{1\cdot 2\ga(\ga + 1)} x^2 + \cdots ;
\]
namely
\[
F \left(-n + l + 1, 1, l + 2, - \frac p q \right).
\]

A hypergeometric series of this type has a continued fraction expansion that is a special case of Gauss's continued fraction for a ratio of two hypergeometric functions (for the latter, see Wall, 1948).  In this particular case, if 
\[
c_k = \frac{(n - k - l)(l + k)}{(l + 2k - 1)(l + 2k)} \frac p q\quad \text{and} \quad  d_k = \frac{k(n+k)}{(l + 2k)(l + 2k +1)} \frac p q ,
\]
then the resulting expansion is
\[
S = \cfrac{1}{  1 - \cfrac{c_1}{  1 +\cfrac{d_1}{  1 - \cfrac{c_2}{  
      \begin{array}{@{}c@{}c@{}c@{}}
        1 + {}\\ &\ddots\\ &&{} -      \cfrac{c_{n-l-1}}{  1    + \cfrac{d_{n-l -1}}{  1     }}  \; .
      \end{array}  
}}}}
 \]
For each $k \geq 1$, let 

\[
C_k = \cfrac{1}
{  1 - \cfrac{c_1}{  1 + \cfrac{d_1}
{      
\begin{array}{@{}c@{}c@{}c@{}}
        1 - {}\\ &\ddots\\ &&{} -      c_k      
\end{array}  
}
}}
\qquad \text{and} \qquad
D_k = \cfrac{1}
{  1 - \cfrac{c_1}{  1 +\cfrac{d_1}
{      
\begin{array}{@{}c@{}c@{}c@{}}
        1 - {}\\ &\ddots\\ &&{} -      \dfrac{c_k}{1 + d_k}      
\end{array}  
}
}}.
\]
We will refer to the $C_k$ and $D_k$ as ``$C$-convergents'' and ``$D$-convergents''.  
It is not hard to see that the Markov method brackets $S$ by a $C$-convergent and a $D$-convergent.  In particular, it is not hard to show  the convergents exhibit a two-step ping-pong type of behavior, satisfying the successive sequence of inequalities
\[
C_2 < D_2 < C_4 < D_4 < \dots < S <  \dots < D_3 < C_3 < D_1 < C_1
\]
(until $k = n - l - 1$, when $D_k = S$).    

Because it is surprisingly hard to find an actual proof of the alternating behavior of the convergents in the literature, we give the details here.  

\begin{proof}
 Because $l > np$, is easy to see that
\[
1 > c_1 > \dots > c_k > \dots  > c_{n-l} = 0.
\]
Let
\[
\om_k := \cfrac{c_k}
{  1 + \cfrac{d_k}
{  1 - \cfrac{c_{k+1}}
{  1 + \ddots}
}
}
\]
Note that by definition one has that
\[
S = \frac{1}{1 - \om_1},
\]
as well as the recursive relation
\[
\om_k = \cfrac{c_k}{1 + \cfrac{d_k}{1 - \om_{k+1} }}.
\]
Furthermore, it is easily seen that $0 < \om_k < c_k$ for $k < n - l$.  

The proof proceeds by exploiting the recursive nature of the $\omega_k$.  Note that for $0 < \al < \beta$,
\[
\cfrac{c_k}{1 + \cfrac{d_k}{1 - \be }} < \cfrac{c_k}{1 + \cfrac{d_k}{1 - \al }}.
\]
Call this the \ti{basic inequality}.  Because $0 < \om_{k+1} < c_{k+1}$, it follows from the basic inequality and the $\om_k$ recursion  that
\[
\cfrac{c_k}{1 + \cfrac{d_k}{1 - c_{k+1} }} < \om_k < \cfrac{c_k}{1 + d_k} < c_k.
\] 
If one then continues this process, thus  extending the continued fraction back by the additional two terms $c_{k-1}$ and $d_{k-1}$, the result is a new set of inequalities with the direction reversed.  (Note this is an immediate consequence of now applying the basic inequality three times.)  Continuing this a total of $k$ times, going all the way to $c_1, d_1$, and then finally passing from bounds for $\om_1$ to bounds for $S = 1/(1 - \om_1)$ (note this last step does not reverse the inequalities), one ends up with
\[
(-1)^k C_{k+1} < S < (-1)^kD_k < (-1)^kC_k.
\]
Furthermore, since
\[
\cfrac{c_{k+1}}{1 + \cfrac{d_{k+1}}{1 - c_{k+2} }} > 0,
\]
once again invoking the basic inequality and arguing as before, one has
\[
 (-1)^{k+2}C_{k+2} < (-1)^kD_k.
\]
It follows the convergents satisfy the sequence of inequalities
\[
C_2 < D_2 < C_4 < D_4 < \dots < S <  \dots < D_3 < C_3 < D_1 < C_1.
\]
\end{proof}

Uspensky illustrates this computational process in detail for the case when $n = 9,000, p = 1/3$, and $l = 3,090$, and starting at $\om_6$ obtains the bounds
\[
0.02161 < P(S_{9000} > 3090) < 0.02175.
\]
As noted earlier, the exact answer to five places  is 0.02170, so in this example the continued fraction approximation is accurate to three places.

\subsubsection{Discussion}

Although Markov's method is at once both elegant and useful, there are some puzzling aspects to Uspensky's presentation.  First, Uspensky computes the convergents of the continued fraction ``from the bottom up''.  This is fine as far as it goes, but inefficient if one then decides to compute further out in the continued fraction because you have to start over from scratch.  

Following the notation in Dudley (1987), consider the interweaved set of convergents $Q_1 = 1$ and 
\[
Q_{2k} = C_k, \qquad Q_{2k+1} = D_k,  \qquad k > 1.
\]
Let $Q_n = A_n/B_n$.  The numerators $A_n$ and denominators $B_n$ of the convergents of an alternating continued fraction satisfy a simple two-step linear recursion which facilitates their computation;  see Dudley (1987, p.\ 589, Theorem 3.1).  Here the formula for the recursion takes the form:  if $X_n$ represents either the numerator $A_n$ or the denominator $B_n$, then for $k > 1$, 
\begin{align*}
X_{2k} &= X_{2k-1}  - \; c_k \, X_{2k-2}, \\
X_{2k+1} &= X_{2k} \quad + \; d_k \; \, X_{2k-1},
\end{align*}
subject to the initial conditions $A_0 = 0, A_1 = 1, B_0 = B_1 = 1$.  One can thus readily compute successive convergents using the recursive formula, knowing $S$ is always bounded by  two successive convergents $D_k, C_{k+1}$ lying  on either side of $S$, stopping when a desired degree of accuracy is reached. 
It is unclear why Uspensky does not mention this, since he was certainly familiar with this phenomenon;  see p.\ 359 of his own book.

For some reason  Uspensky does not give a specific reference for Markov's method.  It in fact  appears in Markov's 1900 book \ti{Ischislenie veroiatnostei} (Calculus of Probabilities), Section 24, pp.\ 150--157;  and Uspensky's treatment follows that of Markov's very closely, including the layout of the illustrative numerical example (although Uspensky uses different values for $n, p$, and $l$).  Markov does note (1900, p.\ 153) the odd $C$-convergents are greater than $S$ and the even ones less, but does not mention their monotonic nature and says nothing about the $D$-convergents.  It is hard to believe however he did not know the full story.  Perhaps, focussing on the computation for a single value of $k$, he did not think it pedagogically desirable to go into such details.  For discussion of Markov's treatment, see Dutka (1981, pp.\ 19--21), Seneta (2013, pp.\ 1105-6).

The use of continued fractions in probability to perform computations goes back to Laplace (Hald, 1998, pp.\ 208--209).  Chebyshev used them to study least squares interpolation (Hald, 1998, pp.\ 525--527), as well as part of a theoretical attack on the moment problem.  Markov's 1884 thesis (``On some applications of algebraic continued fractions''), continued this work of his advisor Chebyshev;  see Uspensky (1937, Appendix II).  Markov also published a short book in French shortly after (Markov, 1888)  in which he used continued fractions to compute tables of the normal integral.

\subsubsection{Later literature}

Markov's book went through four editions (1900, 1908, 1913, 1924), but the section on the method of continued fractions in the last edition (Markov, 1924, Section 18, pp. 104--114), is virtually identical to that  in the 1900 edition.   The second edition is conveniently available in German translation (Liebmann, 1912);  see Section 25, pp. 135--141 for the material on continued fractions. 

Uspensky appears initially to have been virtually the sole source---perhaps the sole source---for discussion of the method in the English literature.  Bahadur (1960) independently came up with a closely related representation of the tail as a product of the lead term and $q$ times a different hypergeometric series,
\[
P(S_n \geq j) = \left[\binom{n}{j} p^j q^{n-j} \right] q F(n+1, 1; j + 1; p).
\]  
(The simple argument invokes a standard relation between the binomial and negative binomial distributions.)  Bahadur then goes on to add:
\begin{quote}
 Another method of using continued fractions to obtain bounds on $B$, which is based [directly on the tail sum] itself rather than [a series deriving using the negative binomial distribution], is given in Uspensky ([2], pp. 52-56).  This method, which is attributed in [2] to Markov, does not appear to be generally known, and might therefore be described here.  
\end{quote}

Much of the subsequent literature on Markov's method is curiously terse when it comes to the monotonic nature of the convergents.   Bahadur says ``it can then be shown'', citing Uspensky, but gives no page reference and, as discussed above, Uspensky is in fact silent on this issue.     Dudley (1987, p. 589), in a very nice paper on ``Some inequalities for continued fractions'', passes over both Bahadur and Uspensky in silence and directly cites Markov (1924, p. 108).  He appears to have felt some frustration with the discussion in the literature about the two-step alternating nature of the convergents, for he says of it that it is ``known, at least in some cases'', citing both Markov (1924, p. 108) and Peizer and Pratt (1968, p. 1452).  He says it ``follows directly from'' the generalized continued fraction representation.  

For a more general discussion of continued fraction expansions for hypergeometric series, see Borwein et al.\ (2005).

\subsection{Bernoulli's theorem}

Nowadays the weak law of large numbers for independent and identically distributed trials is usually proven using Chebyshev's proof.  But this is not the one Uspensky gives for Bernoulli trials.  Instead he says:
\begin{quote}
 Several proofs of this important theorem are known which are shorter and simpler but less natural than Bernoulli's original proof.  It is his remarkable proof that we shall reproduce here in modernized form.  
 \end{quote}
 
\subsubsection{Bernoulli's proof}

The strategy underlying James Bernoulli's original proof is straightforward:  divide the range of the binomial tail probability into blocks, bound the sum of the terms in each block by a corresponding term in a geometric series, and then show the sum of the resulting geometric series satisfies the requisite inequality for all $n$ sufficiently large.  The actual tactical execution of this however requires considerable skill.  Here is an outline of the proof Uspensky gives using his notation.

We seek to show that if there are $m$ heads in $n$ tosses of a $p$-coin, then for any $\ep, \eta > 0$, one has
\[
P\left( \left|\frac m n - p \right| < \ep \right) > 1 - \eta
\]
for all $n$ sufficiently large.

Here is the proof.  Let $\la = \lceil np \rceil, \mu = \lceil np + n\ep \rceil$ denote respectively the ceiling of $np$ and $np + n\ep$ respectively;  that is, the integers $\la, \mu$ such that
\[
\la - 1 < np \leq \la, \qquad  \mu - 1 < np + n \ep \leq \mu.
\]
Let $g := \mu - \la$ (this will be the block size),  let $T_k$ be the probability of $k$ heads in $n$ tosses of a $p$-coin,  and set
\begin{align*}
A &= T_{\la} + \dots + T_{\mu-1},\\
C &= T_{\mu} + \dots + T_n,\\
A_j &= T_{\mu + (j-1)g} + T_{\mu + (j-1)g + 1} + \dots + T_{\mu + jg - 1}, \quad j \geq 1.
\end{align*}
so that
\[
C = A_1 + A_2 + A_3 + \dots .
\]
The proof then breaks down into the following steps:\\

\nin 1.  \tb{Bound the tail by a geometric series}.  For integers $b > a \geq 0$ and $k > 0$ a simple manipulation of inequalities gives
\[
\frac{T_{b+k}}{T_b} < \frac{T_{a+k}}{T_a};
\]
and it is easily seen from this (setting $A_0 := A$) that 
\[
A_j < A_{j-1} \left( \frac{T_{\mu}}{T_{\la}} \right), \qquad j \geq 1,
\]
and therefore
\[
C <  A \left[ \left( \frac{T_{\mu}}{T_{\la}} \right) + \left( \frac{T_{\mu}}{T_{\la}} \right)^2 + \left( \frac{T_{\mu}}{T_{\la}} \right)^3 + \dots \right].
\]

\nin 2.  \tb{Derive a bound for the common ratio of this series}.  For $x \geq \la \geq np$, it is easily seen that $T_{x+1}/T_x < 1$.  Expressing $T_{\mu}/T_{\la}$ as a telescoping product of $g$ terms, retaining only the first $\al \leq g$ of these terms (that is, the length of some sub-block), invoking the first inequality in Part 1 above, and noting that
\[
\frac{T_{\mu - \al +1}}{T_{\mu - \al}} \leq \left(\frac{n - \mu + \al}{n - \al  +1} \frac p q \right),
\]
gives
\[
\frac{T_{\mu}}{T_{\la}} < \left(\frac{n - \mu + \al}{n - \al  +1} \frac p q \right)^{\al}.
\]

\nin 3.  \tb{Use the preceding to show the common ratio of the series is bounded by $\eta$}.  So far $\al \leq g$ can be any sub-block size.  Now choose $\al$ to be the least positive integer such that
\[
\left( \frac{p}{p + \ep} \right)^{\al} \leq \eta.
\]
Further manipulation of inequalities shows that for 
\[
n \geq \frac{\al(1 + \ep) - q}{\ep(p + \ep)}
\]
one has both $g \geq \al$ (so that the last inequality in the previous step holds) and 
\[
\frac{n - \mu + \al}{n - \al  +1} \frac p q  < \frac{p}{p + \ep}.
\]
It then follows immediately that 
\[
n \geq \frac{\al(1 + \ep) - q}{\ep(p + \ep)} \qquad \Longrightarrow \qquad 
\frac{T_{\mu}}{T_{\la}} < \left(\frac{n - \mu + \al}{n - \al  +1} \frac p q \right)^{\al} < \left( \frac{p}{p + \ep} \right)^{\al} \leq \eta.
\]
4.  \tb{Use this bound to show the tail probability $C$ is small than $\eta$}.  Putting all this together then gives (because $A < 1 - C)$
\[
C < A \left( \eta + \eta^2 + \eta^3 + \dots \right) = \frac{A \eta}{1 - \eta}
\quad \Longrightarrow \quad   C < \frac{(1 - C)\eta}{1 - \eta} \quad \Longrightarrow \quad   C < \eta. 
\]
\nin \ti{Quod erat demonstrandum} (as Bernoulli might have said).\\

\nin \ti{Remark}  It is interesting to see Bernoulli uses the device of blocking terms, because this technique was commonly used in the twentieth century to derive limit theorems for sums of independent and identically distributed random variables.\\

In short, Uspensky has provided a clear and complete presentation of an interesting (and impressive) proof by Bernoulli, one that could not be found in any of the textbook literature of the time (at least in English).\footnote{In the 1924 edition of his book, Markov presented Bernoulli's proof but silently included improvements in it due to Nicholas Bernoulli (providing better bounds on the tail probability), as well as a ``modernized'' version of the proof which dispensed with unnecessary restrictions (on $n, p$, and $\ep$).  Strictly speaking, it is Markov's improved and modernized version that Uspensky gives, together with an improvement of his own (a lower bound on $n$ that no longer depended on $p$). Hald (1990, Chapter 16) gives a detailed analysis of Bernoulli's original proof, as well as discussing the contributions of Nicholas Bernoulli and Uspensky.}

\subsubsection{Background}

There is an interesting backstory to Uspensky's presentation of Bernoulli's proof.  The year 1913 was the bicentenary of the appearance of Bernoulli's \ti{Ars conjectandi}, the famous book by Bernoulli in which he gave his proof.  To mark the centenary, Markov arranged a meeting of the St. Petersburg Mathematical Society, as well as commissioning Uspensky (who had an excellent knowledge of Latin) to write a Russian translation of Part 4 of the \ti{Ars};  see Uspensky (1913) and for the meeting itself, Seneta (2013).  Uspensky thus had an intimate knowledge of Bernoulli's book and other material from it also appears scattered throughout Uspensky's book.  (Uspensky in fact owned a copy of the  \ti{Ars conjectandi} which is now part of the Stanford University Library Special Collections.\footnote{But only recently---Steve Stigler was still able to check this out from the main library stacks in 1972!})

Good English translations of the \ti{Ars} are now available;  see particularly Bernoulli (2006), and also Sheynin (2005).  Stigler (1986, pp.\ 63--70) gives an excellent historical account of Bernoulli's theorem and his proof of it.

\subsubsection{Other good material in this chapter}
Since our goal is to focus on the most distinctive and interesting topics in Uspensky, we are not going to go through his book, chapter-by-chapter, in a systematic way.  But to give a sense of the richness of the book, here are some other highlights of just this chapter.  For example, immediately after his treatment of Bernoulli's theorem, Uspensky states and proves \ti{Cantelli's theorem}, something rarely found outside of graduate texts (and often not even then):

\begin{theorem}
For $\ep < 1, \eta < 1$ let $N$ be an integer satisfying the inequality 
\[
N > \frac{2}{\ep^2} \log \frac{4}{\ep^2 \eta} + 2.
\]
The probability that the relative frequencies of an event $E$ will differ from $p$ by less than $\ep$ in the $N$-th
and all the following trials
is greater than $1 - \eta$.
\end{theorem}

\nin Other highlights of the chapter include:
\begin{itemize}
\item A lengthy translation from the \ti{Ars conjectandi}.
\item A survey of experimental verifications of Bernoulli's theorem.
\item A discussion of the Buffon needle problem.\footnote{Including an extension of it to the case of a triple grid;  see Perlman and Wichura (1975, pp.\ 159--162) for a discussion of the statistical aspect of estimating $\pi$ in this case.}
\item The standard Chebyshev proof of the weak law of large numbers given as a three-part exercise.
\item The proof of the Weierstrass approximation theorem,  using Bernoulli's theorem to show the Bernstein polynomials are dense in the continuous functions on a closed and bounded interval $[a, b]$, is sketched in an exercise.
\end{itemize}

\subsection{The Lexis Ratio}

The middle of the 19th century saw an increasing interest in the apparent stability of statistical ratios (say the yearly suicide rates in Paris and Marseilles).  An important issue that arose out of this was whether different sets of frequency data from different populations could be legitimately combined:  that is, whether they represented trials of the same or different phenomena.  Suppose one has $n$ independent series of $s$ independent Bernoulli trials each.  Let $m_1, \dots, m_n$ denote the number of successes in each of the $n$ series of $s$ trials, and let $p$ denote the mean probability of success in all $N := ns$ trials.    In 1876 the German statistician Wilhelm Lexis (1837--1914) defined the \ti{coefficient of dispersion}
\[
Q := \frac{\sum_{i=1}^n (m_i - sp)^2}{Np(1-p)}.
\]
Let $p_{ij}$ be the probability in the $i$-th series that the $j$-th outcome is a success.  Then the mean probability in the $i$-th series is
\[
p_i := \frac{p_{i1} + p_{2i} + \dots + p_{is}}{s}, \qquad \text{and} \qquad p = \frac{p_1 + p_2 + \dots + p_n}{n}.
\]

\nin If $D := E(Q)$, then one can show that 
\[
D = 1 + \frac{s - 1}{np(1-p)} \sum_{i=1}^s (p - p_i)^2 - \frac{1}{Np(1-p)} \sum_{i=1}^n \sum_{j=1}^s (p_i - p_{ij})^2.
\]
There are three natural special cases here:\\

\nin 1.  The \ti{Bernoulli} case:  the $p_{ij}$ have the same value $p$.  In this case $D = 1$ (``normal dispersion'').\\

\nin 2.  The \ti{Lexis} case:  the $p_{ij}$ are constant within strata;  $p_{ij} = p_i$.  In this case the third term in $D$ vanishes but not the second, $D > 1$ (``supernormal dispersion'').\\

\nin 3.  The \ti{Poisson} case:  the $p_{ij}$ vary within a stratum but in the same way from one stratum to another;   $p_{i_1j} = p_{i_2j}$.  In this case the second term in $D$ vanishes but not the third and $D < 1$ (``subnormal dispersion'').\\

There was a substantial literature available in English about this in 1937 (see, for example, Fisher, 1922, Chapter 10;  Forsyth, 1924;   Rietz, 1924, Chapter 6;  Coolidge, 1925, Section 4.2), but Uspensky's treatment introduced contributions by Markov and Chuprov that were unknown in the English literature of the time.\footnote{Julian Lowell Coolidge (1873--1954) was a mathematician at Harvard University who wrote a number of books on mathematics and the history of mathematics.  His 1925 book  was a substantial (if inferior to Uspensky's) book in English at the time on mathematical probability.}

Markov considered an empirical version $\widehat{Q}$ of the dispersion coefficient $Q$ replacing the unknown theoretical probabilities by ones estimated from the data:  if $M = \sum_i m_i$, then
\[
\widehat{Q} = \frac{n(N-1) }{(n-1)} \frac{\sum_{i=1}^n \left( m_i - s \dfrac M N \right)^2}{M(N-M)}.
\]
(If $M = 0$ or $M = N$ then $\widehat{Q} = 1$ by definition.)  Chuprov and Markov were then able to establish the following exact expressions for $E(\widehat{Q})$ and $Var(\widehat{Q})$:
\begin{theorem}
 In the Bernoulli case $E(\widehat{Q}) = 1$ and
 \[
Var(\widehat{Q}) = \frac{2N(N-n)}{(n-1)(N-2)(N-3)}
 \sum_{M=1}^{N-1} \frac{M-1}{M} \cdot \frac{N-M-1}{N-M} \binom{M}{N} p^M (1 - p)^{N-M}.
 \]
\end{theorem}
\nin This immediately gives a simple upper bound for the variance:
\[
Var(\widehat{Q}) <  \frac{2N(N-n)}{(n-1)(N-2)(N-3)};
\]
and, for $n \geq 5$ the even simpler bound
\[
Var(\widehat{Q}) < \frac{2}{n-1}.
\]

For excellent historical accounts of the Lexis ratio and the contributions of Markov and Chuprov, see Heyde and Seneta (1977, Section 3.4),  Ondar (1981), and Stigler (1986, Chapter 6).

\subsubsection{Ships passing in the night}

There was an interesting disconnect between the English and Continental literatures on statistics during the nineteenth and early twentieth centuries.  One simple example is the distribution of $S^2$, the sample variance, in the case of sampling from a normal population:  Helmert had already worked this out in 1876, but this was overlooked in England and independently rediscovered by Student in his  famous 1908 paper on the $t$-statistic (see Zabell, 2008, Section 2.3.1).  Indeed, Helmert's priority was only recognized and acknowledged in the English literature much later (by Karl Pearson in 1931).  

A similar situation held in the case of the Lexis ratio:  although effectively the same as the chi-squared statistic, it was only in 1925 that Fisher observed:
\begin{quote}
It is of interest to note that the measure of dispersion, $\phi^2$, introduced by the German economist Lexis is, if accurately calculated, equivalent to $\chi^2/n$ of our notation.  [Fisher, 1925, p.\ 79]
\end{quote}
Thus, much like Moliere's \ti{bourgeois gentilhomme} who (as the Wikipedia puts it) was ``surprised and delighted to learn that he has been speaking prose all his life without knowing it'', both groups of statisticians had been talking about the same thing without realizing it;  and one could still find papers in \ti{Biometrika} in the 1930s and 1940s computing the exact expectation of the chi-squared statistic, in ignorance of the earlier contributions of Markov and Chuprov.

But--surprisingly--even Uspensky himself was not immune to this breakdown in communication.  After discussing the exact expectation and variance of $\widehat{Q}$, Uspensky went on to write (p.\ 219):
\begin{quote}
 It would be greatly desirable to have a good approximate expression for the probability of either one of the inequalities 
 \[
 \widehat{Q} \geq 1 + \ep \qquad \text{or} \qquad \widehat{Q} \leq 1 - \ep.
 \]
 But this important and difficult problem has not yet been solved.
\end{quote}
This curiously overlooks  not only the work of R.\ A.\ Fisher (with which Uspensky was in fact well acquainted), but also Markov (who had worked out the limiting distribution of $\widehat{Q}$ in 1920).

\subsubsection{Uspensky on error estimation and robustness}

When later discussing the closely related case of Pearson's chi-squared test for the goodness of fit of a vector of observed frequencies to a prescribed vector of multinomial probabilities, Uspensky argued (p.\ 327)  that because ``the degree of approximation [for the test] remains unknown'', the ``lack of information as to the error incurred by using an approximate expression \dots renders the application of this "$\chi^2$-test" devised by Pearson somewhat dubious''.  (This issue was later addressed in part by Cochran, 1952, among others.)  Of course similar concerns can also be raised about the use of $t$-statistics and the other elements of normal sampling theory, and here too Uspensky expressed (p.\ 345) reservations about their use:
\begin{quote}
 The various distributions dealt with in this chapter are undoubtedly of great value when applied to variables which have normal or nearly normal distribution. Whether they are always used legitimately can be doubted. At least the ``onus probandi'' that the ``populations'' which they deal with are even approximately normal rests with the statisticians.
\end{quote}
Such robustness issues had been raised earlier by Egon Pearson in his review of the second (1928) edition of Fisher's \ti{Statistical Methods for Research Workers};  see Pearson (1990, pp.\ 95--101), Zabell (2008, p.\ 4).  Nowadays there is a general consensus that this is less of a concern in the case of estimates or tests for means, but much more so in the case of variances.

In general, Uspensky is typically interested not just with limit theorems as approximations, but also providing an estimate of the magnitude of the error (as in Markov's method of continued fractions).

\subsection{Some other interesting gems}

There are many other interesting topics in Uspensky rarely found in other books, even today.  These include:
\begin{itemize}
\item The use of difference equations in solving problems (Chapter 5).
\item  An analysis of the error in the normal approximation to the binomial (Chapter 7).\footnote{Dutka (1981, p.\ 16)  describes this as a ``very sophisticated version of Laplace's analysis'', and notes the ``considerable technical difficulties'' involved in its proof.  See Seneta (2013, pp.\ 1112--1114) for historical context and discusion.  Interestingly, Littlewood late in life became interested in this problem, although his 1969 paper on it contains several (fixable) errors;  see McKay (1989). 
}
\item  A bound on the error for the Poisson approximation to the binomial (Chapter 7, Exercise 9, pp.\ 135--136).
\item Gambler's ruin with unequal stakes (Chapter 8, pp.\ 143--146).
\item De Moivre's formula for the mean absolute deviation of the binomial \qquad  (Chapter 9, p.\ 176, Exercise 1).\footnote{See Diaconis and Zabell (1991) for the history of this most interesting result.  As noted there, the formula was independently rediscovered several times, including twice after Uspensky's book appeared.}
\item Bernstein's inequality (Chapter 10, pp.\ 204--205, Exercises 12--15).
\item Moment inequalities (Chapter 13, p.\ 278, Exercise 3).
\item A detailed reconstruction of Liapounov's proof of the Liapounov central limit theorem (Chapter 13, pp.\ 284--296).\footnote{In part ahistorical;  see Siegmund-Schultze (2006) for a critique and connections with the later work of von Mises.}
\item A central limit theorem for  two-state Markov chains (Chapter 13, pp.\ 297--302).
\item A crash course in mathematical statistics (Chapter 16).\footnote{Uspensky gives detailed mathematical derivations of the distributions of  $\chi^2, S^2, t, r$, and $\tanh(r)$ for samples from a normal population.  This appears to be have been the first textbook in English to do so.  Although all of these were of course available in the research literature, a similar exposition two years earlier in \ti{The American Mathematical Monthly} (Jackson, 1935) justified itself on the grounds of ``bringing together items that are scattered through a variety of books and journals, and in supplying explanations which in one account or another may have to be read between the lines''.  But even here Uspensky played a role:  in his acknowledgements at the beginning of his paper Jackson wrote:  ``In preparing the paper for publication the author has derived profit from remarks made by Professors Hotelling and Uspensky in the discussion following the oral presentation''.

Jackson was Dunham Jackson (1888--1946), a mathematician at the University of Minnesota known primarily for his work in approximation theory.  He was a President of the MAA, active in the AMS, a Fellow of the American Physical Society, and a member of the National Academy of Sciences;  see Hart (1948 and 1959).}
\item Stirling's formula with bounds and the evaluation of definite integrals (Appendix 1).\footnote{The log convexity argument in Section 3 may  go back to Bohr-Mollerup.} 
\item The method of moments and its applications (Appendix 2).\footnote{This is a difficult topic, with very few expository presentations for non-experts.  We do not know of anything comparable.}
\item Kuzmin's solution to Gauss's continued fraction problem (Appendix 3).
\end{itemize}

Several of these topics warrant brief mention because they crop up in later literature.

\subsubsection{The problem of runs}

One illustration of the use of difference equations Uspensky gives is to solve the \ti{problems of runs} (Sections 3--8,  pp.\ 77--84):  to find the probability $y_n$ of at least $r$ consecutive successes in $n$ tosses of a $p$-coin.  The $y_n$ are easily seen to satisfy the recursion
\[
y_{n+1} = y_n + (1 - y_{n-r})qp^r.
\]
If we let $z_n := 1 - y_n$, this gives us the recursion
\[
z_{n+1} - z_n + qp^rz_{n-r} = 0,
\]
and this in turn enables us to find the generating function of the $z_n$.  Uspensky (pp.\ 78--79) gives the detailed argument, finding the generating function of the $z_n$ to be
\[
\phi(\xi) = \frac{1- p^r\xi^r}{1 - \xi + qp^r\xi^{r+1}},
\]
which can be expressed as a power series in $\xi$ ``according to the known rules''.  The $z_n$ are the coefficients of this power series.  

Of course for this to be useful one want a formula for the coefficients.  If one sets
\[
\be_{n, r} :=  \sum_{k = 0}^{\lfloor{\frac{n}{r+1}}\ \rfloor} (-1)^k \binom{k}{n - kr} (qp^r)^k, 
\]
Uspensky  says one can show ``without any difficulty'' that
\[
z_n = \beta_{n, r} - p^r \beta_{n-r, r}. 
\]
For the details of this last   calculation, see Ethier (2010, p.\ 41).  Although the generating function for the $z_n$ was known long before, going back to Laplace (see Todhunter, 1865, pp.\ 539--541;  Hald , 1990, pp.\ 417--421). Uspensky's formula for its coefficients may originate with him.  

Obviously this formula for $z_n$ can be  challenging to compute for large $n$, and Uspensky devotes several pages (79--84)  giving approximate formulas for $z_n$ when $n$ is large, drawing on the results in an earlier (but uncited) paper of his (Uspensky, 1932).

\paragraph{Prehistory}

The problem of runs was first solved by De Moivre (1738, Problem 88, pp.\ 243--248;  1756, Problem 74, pp.\ 254--259). De Moivre gave the following algorithm to compute the probability of a run:  if $q = 1-p$ and $c = p/q$, then expanding 
\[
\frac{p^r}{1 - q - cq^2 - c^2q^3 \dots - c^{r-1}q^r}
\] 
in powers of $q$ and taking the sum of the first $n-r+1$ terms ``will express the probability required''.\footnote{As Todhunter notes, De Moivre gives an erroneous value for $c$ in his algorithm, but his examples are correctly computed.}  De Moivre did not explain where his formula came from, but Todhunter (1865, Section 325, pp.\ 184--6) helpfully sketched a derivation.  

\subsubsection{Gambler's ruin with unequal stakes}

Uspensky's Chapter 8 (``Further Considerations on Games of Chance'') begins with a discussion of the standard topic of computing the absorption probabilities for gambler's ruin with equal stakes, but then (p.\ 143) turns to the nonstandard question of what happens when the stakes are unequal: 

\begin{quote}
Two players $A$ and $B$ play a series of games, the probability of winning a single game being $p$ for $A$ and $q$ for $B$, and each game ends with a loss for one of them. \dots  [If] the stakes for $A$ and $B$ measured in a convenient unit are $\al$ and $\be$ and their respective fortunes are $a$ and $b$, find the probabilities for $A$ and $B$ to be ruined in the sense that at a certain stage the capital of $A$ become less than $\al$ or that of $B$ less than $\be$.
\end{quote}

If $y_a$ is the probability of ruin for $A$, Uspensky (pp.\ 143--146) shows in the case of a fair game (i.\ e.\, $p\be = q\al$) that $y_a$ satisfies the inequalities
\[
\frac{b - \be + 1}{a + b - \be + 1} \leq y_a \leq \frac{b}{a + b - \al + 1},
\]
as well as a slightly more complicated set of inequalities in the unfair case.

Here, once again, Uspensky is drawing on ``an ingenious method developed by Markov'' (Markov, 1912, pp.\ 142--146).  The problem had been studied earlier by Rouch\'{e} (1888) and discussed in Bertrand's 1888 textbook, the crucial expression being the equation
\[
pz^{\al + \be} - z^{\al} + q = 0.
\]
Markov, in an obscure paper in 1903, showed that all the roots to this equation contribute to the probability of ruin and derived the inequalities Uspensky gives.  Markov's method also appears briefly in Feller (1957, p.\ 366 of the 1968 edition), 

Since Uspensky the problem has continued to generate a modest literature.  Feller generalized Uspensky's result in 1950  (see  p.\ 366 of the 1968 edition of his book).   Hillary Seal (1966) traces the history on this up to 1957.  Subsequent literature includes Ethier and Khoshenevisan (2002) and Gilliland et al. (2007), the later giving exact formulas in terms of the roots of $pz^{\al + \be} - z^{\al} + q = 0$.  Ethier (2010, Sections 7.2-3) discusses these generalizations at length and reviews their history (pp. 273--4).

\subsubsection{Bernstein's inequality}

Let $X_1, X_2, \dots, X_n$ be independent random variables,,  such that $EX_j = 0$ and $\si_j^2 := Var X_j < \infty$ for $j = 1, \dots, n$.  \ti{Bernstein's inequality} states that if  $S_n := X_1 + \dots + X_n, B_n^2 = Var S_n$, and for some $c > 0$ one has
\[
E\left[\left|X_j^k \right| \right]\leq k! \, \frac{\si_j^2}{2} \, c^{k-2}
\]
for $j = 1, . \dots , n$ and $k > 2$, then
\[
P\left(\left|S_n \right| > t \right) < 2 \exp\left(- \frac{t^2}{2B_n + 2ct} \right). 
\]
In particular, if the $X_j$ are uniformly bounded, that is, $\left|X_j \right| \leq M$ for some $M > 0$, then one can take $c = M/3$.  This is (like a number of other  interesting results in Uspensky) established via a sequence of exercises;  see Uspensky (1937, pp.\  204--205, Problems 12--15).  Nowadays better results are available, but Bernstein's inequality remains a simple and useful upper bound.

The inequality was proved by \ti{Sergei Natanovich Bernstein} (1880--1968), and appears in  his 1927 book \ti{Theory of Probability} (as well as an earlier 1924 paper), but even in 1962 these were 
hard to find in the US. (Bennett, 1962, 
listing previous references to the inequality in the English literature, 
described  them as being ``unobtainable'').  Uspensky's account in his book---with one exception---remained the sole English source for the inequality for many years.  Bennett (1962, p.\ 35), discussing the inequality a quarter of a century later, could only find six previous instances where it was mentioned in the English literature;  only one of these, Craig (1933), predates Uspensky's book or gives a derivation.  But this turns out to be the exception that proves the rule.  In his paper on the inequality Craig says (p. 94):

\begin{quote}
Another interesting and important attempt in this direction [that is, Chebychev's inequality for a sum of independent random variables] due to
 S.\ Bernstein seems to have generally escaped attention in the English-speaking world, at least, since it has been published only in Russian.
\end{quote}
How then did Craig learn of it?  He tells us in a footnote:  ``The present account of this work of Bernstein is taken from a lecture of Professor J.\ V.\ Uspensky''\,!\\

For other citations of Bernstein via Uspensky, see Blackwell (1954, p.\ 397) and Bahadur (1966, pp.\ 578--579).\footnote{In a note added in proof, Blackwell reports Ted Harris had drawn his attention to the Bernstein inequality, cites Uspensky, and notes that using it would  yield his geometric rate of decay result under weaker conditions and with a slightly better rate.  Bahadur told a colleague that he would not have been able to write his 1966 paper if he had not known of this result in Uspensky.}

\section{Uspensky's two other textbooks}

Uspensky wrote two other textbooks during his time at Stanford, one on number theory (Uspensky and Heaslet, 1939) and one on the theory of equations (Uspensky, 1948).  While our primary focus has been on Uspensky's book on mathematical probability, his other two textbooks are not without interest and so we pause to briefly comment on them.

\subsection{Elementary Number Theory}

Uspensky was first and foremost a number theorist (an interest he shared with his advisor Markov).  So not surprisingly his book  on number theory also contains many interesting examples and topics.  Here are a few:

\begin{itemize}
\item Lam\'{e}'s theorem (the number of divisions needed to find the gcd of two numbers is at most five times the number of digits in the smaller number, pp.\ 43--45).
\item Bonse's inequality (if $p_n$ is the $n$-th prime, then $p_{n+1}^2 < p_1 p_2 \dots p_n , \, n \geq 4$, pp. 86--89).
\item Meissel's formula for computing the number of primes less than a given number (pp.\ 120--124).
\item Kummer's proof of the law of quadratic reciprocity (pp.\ 375--379).
\item Dickson's proof of the four-square theorem (pp.\ 379--386).
\item Chapters on the Bernoulli numbers, including the Voronoi congruences (Chapter 9) and Liouville's methods for deriving properties of arithmetical functions, culminating in an elementary proof of Gauss's theorem characterizing integers  expressible as a sum of three squares (Chapter 13).
\item Appendices on magic squares (pp\ 159--172), calendar problems (pp.\ 206--221), and card shuffling (pp.\ 244-248).
\end{itemize}

\nin We cannot resist briefly discussing this last topic.

\subsubsection{Card shuffling}

Uspensky and Heaslet (1939, pp.\ 244--248) is in effect a standalone article on shuffling cards.  They treat perfect shuffles and the ``Monge'' or ``over-under'' shuffle.  Consider a deck of $2n$ cards (for example, $2n = 52$) in order $1, 2, \dots, 2n$ from the top down.  The deck is cut exactly in half and then the two halves riffle shuffled together as shown below:
\begin{align*}
&1 --& && &1 --\quad-- 5& && &\qquad-- 5& && &5 -- \\
&2 --& && &2 --\quad-- 6& && &1 -- &     && &1 -- \\
&3 --& && &3 --\quad-- 7&  && &\qquad-- 6& && &6 -- \\
&4 --& &\rightarrow& &4 --\quad-- 8& &\rightarrow& &2 -- &  &\rightarrow&      &2 -- \\
&5 --& && && && &\qquad--7& && &7 -- \\
&6 --& && && && &3--& && &3 -- \\
&7 --& && && && &\qquad--8& && &8-- \\
&6 --& && && && &4--& && &4 -- \\
\end{align*}

\nin This is a perfect ``in shuffle'' practiced by gamblers and magazines since (at least) 1743.

The question is, how many perfect shuffles ($r$) are required to recycle the cards---that is, to bring them back to their original order?  When one of us was quite young (PD), we figured this out ``the long way'', by actually shuffling until the cards returned to order.  Here is some of the data:

\bigskip
\begin{table}[htbp]
  \centering
  \begin{tabular}{@{} c | ccccccccccccc @{}}
    $2n =$ & 2 & 4 & 6 & 8 & 10 & 12 & 14 & 16 & 18 & 20 & 22 & 24 & 26\\ 
\midrule
    \quad $r =$ & 2 & 4 & 3 & 6 & 10 & 12 & 4 & 8 & 18 & 6 & 11 & 20 & 18\\ 
  \end{tabular}

\end{table}

\begin{table}[htbp]
  \centering
  \begin{tabular}{@{} c | ccccccccccccc @{}}
    $2n =$ & 28 & 30 & 32 & 34 & 36 & 38 & 40 & 42 & 44 & 46 & 48 & 50 & 52\\ 
\midrule
    \quad $r = $ & 28 & 5 & 10 & 12 & 36 & 12 & 20 & 14 & 12 & 23 & 21 & 8 & 52\\ 
  \end{tabular}
  \caption{Number of shuffles $r$ required for a deck of size $2n$}
  \label{tab:label}
\end{table}

In their book Uspensky and Heaslet prove that after a single shuffle card $i$ moves to position $2i \pmod{2n+1}$, so after $k$ shuffles card $i$ moves to position $2^k i \pmod{2n+1}$.  It follows the deck returns to its original order the first time $2^k i \equiv 1 \pmod{2n+1}$.  Suppose $2n+1$ is prime (for example, when $2n = 52$).  Observe that in this case $2n$ shuffles are  \ti{sometimes} required before the deck returns to its original order.  Does this happen infinitely often?  Nobody knows---this is a special case of the \ti{Artin conjecture} (that 2 is a primitive root for infinitely many $ p$).  The conjecture is known to be true on the generalized Riemann hypothesis but this is a million dollar problem.

As far as we know, Uspensky and Heaslet were the first both to discover and prove that the number of shuffles required to return the deck to its original configuration is the order of two.  For more on the history of perfect shuffles and their applications, see Diaconis, Graham, and Kantor (1983), Diaconis and Graham (2012).  Uspensky and Heaslet's Appendix (pp.\ 245--248) also treats the ``Monge shuffle'' with similar results.  This particular shuffle was first analyzed by the French mathematician Gaspard Monge in 1783.\\

Chapter 13 (``Liouville's Methods'') was a bit of an indulgence on Uspensky's part, because it discussed a major research interest of his, the results being ``elementary'' only in terms of the methods used.  The chapter derived a variety of very general arithmetical identities using basic number theory arguments rather than appealing to the theory of elliptical functions.  Uspensky then used special cases of these identities to derive the Gauss and Jacobi theorems on the representation of an integer as a sum of three or four squares.

 Of  course Uspensky's book did not touch on most of his research (it was, after all, on elementary number theory), but before leaving it, we discuss one other number theoretic result of Uspensky's, one not treated in his book but which remains surprisingly fresh even today.
\subsubsection{A game}

For $\al$ a positive real number, define the \ti{spectrum} of $\al$ to be the sequence of integers
\[
\al_n := \{\lfloor n \al \rfloor : n = 1, 2, 3, \dots \},
\]
where $\lfloor x \rfloor $ is the floor of  $x$.   Such sequences are termed \ti{Beatty sequences} (named after Samuel Beatty, 1881--1970).  A classical theorem attributed to Beatty says that if $\al, \be >1$ are a pair of irrational numbers such that $1/\al + 1/\be = 1$, then the spectra of $\al$ and $\be$ are disjoint and cover the natural numbers.  (This appears on p.\ 98 of Uspensky and Heaslet as Exercise 9!)  It is natural to ask about a similar triple $\al, \be, \ga$ (or more).  Uspensky (1927) proved coverage with disjoint spectra only occurs in the case of  two spectra.  His proof is involved, using Kronecker's approximation theorem.  An elementary proof due to Ron Graham (1963) led to fascinating further theory and many still open problems;  see Graham et al. (1978), Graham and O'Bryant (2005).  \\

Curiously, even before Beatty the result was stated by John William Strutt (Lord Rayleigh).  Rayleigh (Strutt, 1894, p. 123) says:
\begin{quote}
Thus, if $x$ be an incommensurable number less than unity, one of the series of quantities $m/x, m/(1 - x)$, where
$m$ is a whole number, can be found which shall lie between any given consecutive integers, and but one such quantity can be found.
\end{quote}

Rayleigh did not provide a proof.  Beatty (1926) posed it as a problem in the \ti{American Mathematical Monthly}, and several solutions were then submitted by a number of contributors (Beatty, 1927).  Apparently independently, Willem Abraham Wythoff (1865--1939) proposed in 1907 a modification of the game Nim in which the ``cold'' (losing) positions are parametrized by the pair of complementary Beatty sequences generated by $\phi$ and $\phi^2$, where $\phi$ is the golden ratio $(1 + \sqrt{5})/2$.  (Curiously, Uspensky discusses Nim at some length in his book on number theory, pp.\ 16--19, but does not note this connection.)

\subsection{The theory of equations}

Uspensky's last book, \ti{Theory of Equations}, appeared posthumously.  Sent to the publisher in December 1946, only a month before he died (on January 27, 1947), the manuscript was seen through the press by Uspensky's former students Max Heaslet and Carl Olds, and appeared in 1948.  It is less distinctive than the other books:  it was, as Uspensky notes in the preface, ``elementary in nature and, with few exceptions, contains only material customarily included in texts of this kind''.  It was longer than other texts on the subject then currently in use, being designed for self-study if desired.

But---Uspensky being Uspensky---there were still elements of novelty, ones in which ``the exposition differs considerably from custom''.  In Chapter 1, on complex numbers, ``the superficial approach so common in many books'' was replaced by the rigorous definition of a complex number as an ordered pair $(a, b)$ of real numbers.  In Chapter 6 a method for separating real roots was given based on Vincent's Theorem.\footnote{Vincent's theorem, named after Alexandre Vincent (1797--1868) and published by him in 1834, appears to have been entirely forgotten until  Uspensky's discussion of it in his book  (see Chapter 6, Section 12 and Appendix 2).    Uspensky describes the method as ``very efficient'' and ``much superior in practice to that based on Sturm's Theorem'', adding  he ``believes that no other book mentions this method, which he invented many years ago and has been teaching to his students for a number of years'' (p.\ v).  Uspensky's Appendix 2 gave a sharper version of Vincent's result, providing a bound on  the number of steps in the algorithm.   
Shortly after, Ostrowski (1950) was able to improve Uspensky's bound (in a paper which was itself however overlooked for many years).  
For later literature, see Krandick and Mehlhorn (2006).}     
In Chapter 7, on the approximate evaluation of roots that had been separated, Sections 2--7 were devoted to the original form of Horner's method, ``which unfortunately has disappeared from American texts''.  In Chapter 9 determinants were introduced using Weierstrass's approach, based on their ``characteristic properties'' rather than ``formal definition''.  (This is similar to Artin's axiomatic development of them, which can be found in some editions of Lang's calculus and linear algebra textbooks.)  And in an appendix Uspensky gave Gauss's fourth proof of the fundamental theorem of algebra.

\begin{center}
***
\end{center}

All of Uspensky's books are marked by exemplary intellectual sharpness:  the presentation is easy to follow, the proofs are complete, readable, and coherent, in some cases results are given in greater generality than is ordinarily the case.  The writing, thinking through of the material, order of presentation, and choice of examples is both engaging and natural.  The references to prior literature are invariably to the masters, not the pupils;  books and papers are not just cited for show but are ones Uspensky has both read and mastered.

Uspensky grew up in the classical Russian tradition;  and---perhaps as a result---he exhibits an impressively broad mathematical culture.  This suggests some natural questions:  who was he, what was his mathematical background, how did he come to the US?

\section{A closer look at Uspensky}

Detailed information about Uspensky is not easy to find in English, so in this section we discuss his life, based in part on family records, unpublished material in the Stanford archives, material available only in Russian, and many passing references to him throughout the English literature.

\subsection{Uspensky in Russia}

\ti{Yakov Viktorovich Uspenskii} was born on April 29, 1883, the fourth of five children, in Urga, Mongolia.\footnote{Uspensky almost invariably signed his  books and papers in English as ``J.\ V.\ Uspensky''.  The sole exception is his \ti{Annals of Mathematics} paper of 1927, in which he uses the name ``James V.\ Uspensky''.  
His Declaration of Intention for naturalization, dated March 24, 1930, gives his full name as ``James Viktorovitch Uspensky'';  but his 1942 Registration Card for the draft shortened this to ``James Victor Uspensky''.
}   (His father Viktor Matveevich Uspenskii,1845--1901,---a career diplomat---was the Russian consul there.\footnote{See  ru.wikipedia.org/wiki/\foreignlanguage{russian}{Uspenski\u{i},\_Viktor\_Matveevich} (accessed on March 8, 2020) for further information about Uspensky's father.})
When Uspensky was seven his mother and her children moved to St.\ Petersburg to enable the children to receive a proper Western education, while his father (then in the Chinese province of Sinkiang where he was Russian consul-general for Western China) remained at his diplomatic post except when on vacation in Russia.  

Uspensky attended a classical \ti{gymnasium} in St.\ Petersburg;  it was here he learned both Greek and Latin, knowledge which was to prove useful later in his historical studies.  Although initially interested in astronomy, early on he became interested in mathematics;  by the time he graduated he had an excellent grounding in the differential and integral calculus, and was able to read books on astronomy and theoretical physics.  He graduated from the \ti{gymnasium} in 1902 with distinction, receiving a gold medal in recognition of his achievements.

\subsubsection{The Russian university system}

In order to understand Uspensky's university career, a few words about the Russian system then in place may be helpful.  Starting at the beginning of the 19th century, there were three academic degrees:  the \ti{kandidat, magister}, and \ti{doktorat}.  The \ti{kandidat} (later the \ti{diplom}) was a bachelor's degree, awarded after completing one's studies and passing a set of examinations;  the \ti{magister} and \ti{doktorat} were graduate degrees, awarded after writing and publishing a thesis, which then had to be defended in a public oral examination.  (In the case of the \ti{magister} one also had to take and pass a set of examinations prior to the writing of the thesis.)  See generally Sanders (1993).

The \ti{magister} (Latin for ``master'' or ``teacher'', as in \ti{Magister Ludi}) was the equivalent of today's PhD, and required in order to teach as a \ti{privat dozent}, that is, an instructor who took private pay students.\footnote{The \ti{magister} is sometimes referred to as a master's degree, but this is misleading given the master's current  academic status in relation to the PhD.  In contrast, in the medieval universities of Europe the titles Master and Doctor were effectively synonymous, and which degree was conferred depended on on both the university and the faculty within it (arts, law, theology,  medicine);  see Verger (2003, p.\ 146), Rashdall (1895, pp. 21--22).} 
 The \ti{doktorat} was a second doctoral decree, required the writing of a second thesis, and was usually necessary in order to become a professor at a university.  There were finally in turn two grades of professorship, termed ``extraordinary'' and ``ordinary'', roughly equivalent to being an Associate and Full Professor, respectively.\footnote{The terminology derives from the Latin words \ti{extraordinarius } and \ti{ordinarius}, where the ``extra'' connotes difference  (as in extralegal and extraterrestrial) rather than superiority.  The system was similar to the German one, where the two ranks were \ti{au\ss erordentlicher} and \ti{ordentlicher Professor}, professorships without and with a chair, respectively.}

\subsubsection{Uspensky's university career}

Uspensky spent his entire academic career in Russia at the University in what was then St.\ Petersburg.\footnote{The university underwent several name changes during the period we will be discussing:  Saint Petersburg Imperial University (1821--1914), Petrograd Imperial University (1914--18), Petrograd State University (1918--24), and Leningrad State University (1924--91).  The city of St.\ Petersburg was renamed Petrograd immediately after the outbreak of war because both the ``Sankt" and the ``burg'' in ``Sankt-Peterburg'' were German words.  The second name change followed the overthrow of the Czarist regime in 1917, and the third reflected the change of the name of the city from Petrograd to Leningrad five days after Lenin died on January 21, 1924.}   He was an undergraduate there from 1903--1906, years which overlapped the first Russian revolution of 1905 and saw a considerable disruption of Russian academic life.  But unlike many of his fellow students he did not become involved in politics and devoted himself entirely to his studies.  (His apolitical nature may in part explain his escaping relatively unscathed during the upheavals of the 1920s.)  Even at this early stage his talent was evident:  while still an undergraduate he wrote his first paper, the first rigorous proof that the cyclotomic ring $\mathbb{Z}[\zeta_5]$ ($\zeta_5$ a primitive 5-th root of unity) was Euclidean (Uspensky 1906);  this work was regarded as sufficiently important that a reworked version of it written in French was published in the \ti{Mathematische Annalen} in 1909.  In 1906 Uspensky graduated, receiving a \ti{diplom} of the first degree (that is, with honors), was awarded a scholarship, and immediately began his graduate studies.  By 1908 he completed his course work, passed his examinations, and began working on his thesis under the direction of Markov.  The thesis (``Some applications of continuous parameters to number theory'') was completed, approved, and published in 1910;  after passing his oral examination the next year in 1911, Uspensky was then formally awarded the degree of \ti{Magister of Pure Mathematics}.\footnote{Most sources say Uspensky received his PhD or doctorate in 1910;  although essentially correct, the actual title and year the degree was awarded are as indicated.}

With his \ti{magister} in hand, Uspensky became a privat-dozent at St.\ Petersburg in 1912 (although he also supplemented his income during this period by teaching at other institutions:  the Institute of Railway Engineers, 1907--29, the Higher Women's Courses, 1911--17).  In 1915 he became an extraordinary professor at Petrograd Imperial University;  this was unusual because he had not yet been a awarded his \ti{doktorat}, but presumably reflected both his ability and a wartime shortage of personnel.  He was in turn rapidly promoted to ordinary professor in 1917;  this again reflected his ability but  the abolition of the \ti{doktorat} around this time may also have played some role as well.  

In 1921 Uspensky was elected to the Russian Academy of Sciences (as a replacement for Liapunov, who had died in 1918), his election supported by Markov, Steklov, and Krylov.  Their report recommending Uspensky's election (Markov et al., 1921), besides giving a complete bibliography up to that time, cited several of his papers as justifying this honor, in particular a paper (Uspensky, 1920) in which Uspensky had derived (independently of Hardy and Ramanujan, 1918) the asymptotic formula for the number of partitions of an integer,
\[
p(n) \sim \frac{1}{4n\sqrt{3}} \exp \left( \pi \sqrt{\frac 2 3 n} \right),
\]
as well as an estimate for the error.\footnote{Although Uspensky did not derive the full divergent series expansion for $p(n)$ that Hardy and Ramanujan found using their celebrated ``circle method'', he arrived at the simpler formula displayed above by deriving the first term in the expansion together with an estimate of the error, namely
\[
p(n) = \frac{ e^{ \pi \sqrt{ \frac{2}{3} \left( n - \frac{1}{24} \right) } } } {4 \sqrt{3} \left( n - \frac{1}{24} \right) }
\left( 1 - \frac{\sqrt{3}}{\pi \sqrt{2n - \frac{1}{12}}} \right) + \rho_n e^{\frac{\pi}{\sqrt{6}} \sqrt{n} },
\]
where the $\rho_n$ are bounded (Uspensky, 1920, p.\ 209, but correcting a typographical error in Uspensky's formula where a $\sqrt{\pi}$ appears instead of a $\pi$).  This corresponds to Hardy and Ramanujan's equation (1.55) on p.\ 82.  As they noted, this ``is an asymptotic formula of a type far more precise'' than the simpler one, and it was ``with considerable surprise that we found what exceedingly good results the formula gives for fairly large values of $n$''.  (Uspensky's bound on the error is less precise than Hardy and Ramanujan's, however, since in fact $\rho_n = O(1/n)$.  See generally DeSalvo, 2021.)}
(Because of wartime and revolutionary conditions the Hardy and Ramanujan paper was unknown in Russia, even in 1920.)

\subsection{From Petrograd to Palo Alto}

Thus by 1921 Uspensky was a Professor at the Petrograd State University and a member of the Russian Academy of Sciences, an established and senior member of the Russian mathematical community.  And yet eight years later he chose to emigrate to the US.  It is natural to suppose that this was connected with the Russian Revolution---and this was indeed the case---but why did it take so long?  There are some complexities here.

In 1918 \ti{Sovnarkom} (the Council of People's Commissars, that is, the newly established Bolshevik government) began to pass decrees designed to ensure new educational opportunities for workers and peasants.  These measures were resisted by some university faculty concerned that unqualified students might be admitted.  Uspensky for example wrote that
\begin{quote}
Taking into account the fact that to succeed at the university, a student should be adequately trained, prospective students must be admitted at the university in virtue of their knowledge, not their class affiliation or political commitment.\footnote{Central State Historical Archives of St.\ Petersburg. F.\ 7240, Schedule 14, No. 16, L.\ 185 recto;  cited in Nazarov and Sinkevich, 2018, p.\ 6.}
\end{quote}
The response of the Bolshevik government was predictable:  measures were swiftly enacted to ``re-educate the bourgeois professors'' as well as punitive measures such as preventive detention, exile, and even execution (Nazarov and Sinkevich, 2018, p.\ 6).  

These measures were not confined to university professors.  In the summer of 1922, a list of 217 perceived anti-Soviet intellectuals (writers, professors, scientists, etc.) was drawn up and the GPU (\ti{Gosudarstvennoye politicheskoye upravlenie}, the State Political Directorate) proceeded to arrest them beginning on the night of August 16--17.\footnote{The GPU was the  successor of the \ti{Cheka} (the first Soviet secret police, headed by the sinister Felix Edmundovich Dzerzhinsky, 1877--1926) and combined both internal security and foreign intelligence functions.}  Many of those detained were then exiled;  at least 160 of these by ship from Petrograd to Stettin (today the Polish city of Szczecin) on the German ships \ti{Oberb\"{u}rgermeister Hacken} on September 2 and \ti{Preussen} on November 15.  Despite his initially voicing opposition to Soviet university reorganization, Uspensky was not included among these, perhaps in part because he was perceived as being apolitical, perhaps in part because of his abstruse field (as opposed to history, literature, politics, and philosophy), and perhaps because of the prestige accruing from his membership in the Russian Academy of Sciences.  But a chilling message had been sent as to the possible limits of dissent.\footnote{For further information about the 1922 expulsions, see Finkel (2003), Chamberlain (2006).}

In 1924 Uspensky attended the International Congress of Mathematics in Toronto and spoke there about work with his student Boris Venkov (Uspensky and Venkov, 1928);  This may have been Uspensky's first trip outside Russia, and he took advantage of the opportunity by afterwards visiting the University of Chicago and  the University of Michigan, lecturing on Russian contributions to number theory before returning home.  Either because he felt after this trip that publishing his work almost exclusively in Russian had limited its impact in the West, or because he had already begun to contemplate emigrating to the US, from this point on there was a marked shift in how and where Uspensky chose to publish his work.  Although he continued to submit papers to the \ti{Bulletin of the  Academy of Sciences  of the USSR} during the next two years (1925--1926), all but one of the nine---a memorial notice about Steklov primarily of interest to his Russian readers---now appeared in French.\footnote{The Russian Academy of Sciences became the Academy of Sciences of the USSR in 1925, and the name of its \ti{Bulletin} changed accordingly.}  

Uspensky seems to have liked what he saw in the United States, because he returned two years later, teaching at Carleton College in Minnesota during the academic year 1926--1927, followed by short lecture courses at Berkeley and Stanford.  His time at Carleton marks the point at which he first started to publish in English and in US journals (fifteen papers between 1927 and 1935).
During this trip he also met his future wife, Lucile Zander, who was working in the Carleton publicity department;  they married on October 13, 1927, and left for Russia shortly after.

\paragraph{A chat with the OGPU}

When he returned to what was now Leningrad, Uspensky was interrogated by the OGPU (the Joint State Political Directorate, the successor to the GPU after the formation of the USSR in late 1922).\footnote{Royden, who reports this incident, identifies Uspensky's interrogator as a member of the NKVD (the People's Commisariat for Internal Affairs), saying it was a predecessor of the KGB, but this is an anachronism:  the NKVD only acquired the secret police and foreign intelligence functions of the OGPU in 1934.  (In general, the study of the Soviet security apparatus is complicated by its numerous reorganizations and consequent name changes that took place during the period from 1917 to 1995:  the Cheka, GPU, OGPU, NKVD, NKGB, MVD, MGB, KGB, FSB, and SVR.)}  

\begin{quote}
[He was asked] how he liked America.  Uspensky disarmed his interviewer by saying, ``I loved it.  It is a place of great opportunity, and if only I were a young man I would emigrate.  But I am a member of the Academy of Science, and my career is established here.  I am too old to start over again.'' The NKVD  agent evidently reported that Uspensky was reliable and sound in his views.  Thus, when Uspensky did decide to come to America a few years later, he came in style on a Soviet ship with his passage paid for by the government.  
[Royden, p.\ 243]
\end{quote}
(Uspensky's response would have been most ill-advised a decade later at the height of the Great Terror.)

But given his entirely plausible answer to the OGPU---for Uspensky certainly was an established figure in Russian mathematics---why did he leave Russia just a year and a half later?  There are said to have been a number of contributing factors.  The Russian historian Natalia Ermolaeva (1997) reports:

\begin{quote}
Returning to the USSR, Uspensky resumed his numerous duties. However, in the summer of 1929  when he again went on a business trip to the USA, he did not return to his homeland. The decision to emigrate was caused by various reasons. One of them was that during his 2nd trip to America, U. married, and his wife categorically refused to live in the USSR. At the same time, the situation in the mathematical life of the country deteriorated sharply---there was an intensive introduction of Marxism into mathematics accompanied by persecution of scientists, including Nikolai Maximovich G\"{u}nter [a close colleague of Uspensky].
\end{quote}

The final straw appears to have been a tumultuous meeting of the Leningrad Mathematical Society (LMO):  G. G. Lorentz (who was an eyewitness) reports that

\begin{quote}
Arriving at the lecture  I found the entire mathematical Leningrad present.  Uspensky, on the podium, was pointed out to me.   After the lecture, G\"{u}nter, presiding, invited the audience to pose questions. Suddenly Leifert [a Bolshevik, himself later a victim of the Great Terror] climbed on the podium yelling insults at the LMO and G\"{u}nter.  Many students applauded Leifert and shouted. The meeting was dissolved. In the wake of this disastrous event, using their still valid visas, Uspensky and his wife left for the United States.\footnote{Lorentz goes on to add:  ``where he accepted a Stanford University professorship offered to him by Szeg\"{o}''.  Here Lorentz errs:   Szeg\"{o} only left Germany for the US in 1934 and only arrived at Stanford in 1938.  (Szeg\"{o} may of course have played some role from afar in facilitating Uspensky's appointment.)}
\end{quote}

Or maybe not;  perhaps Uspensky had in fact already made up his mind.  An obituary for Uspensky in \ti{The Sheboygan Press} (in Wisconsin, where Uspensky's father-in-law Otto Zander had been the Editor), states:
\begin{quote}
When it came time for Mr. Uspensky to return to Russia  [after his year at Carleton], with his exchange period expired, his wife returned there with him and they lived for a year in Russia, during which time he \ti{completed} arrangements to return to America, leaving his native land to take up citizenship in the new country he had adopted.  [February 3, 1947, p.\ 8;  emphasis added.]
\end{quote}
If this is accurate---and inasmuch as the information in the obituary clearly came from the family there is no reason to doubt it---Uspensky had already decided to emigrate to the US 
\ti{before} he returned to Russia.  (Of course the deteriorating conditions he encountered there on his return, as Stalin increased his grip on the country, would not have helped.)  Indeed, one can be forgiven for conjecturing that Uspensky had even earlier still, when he arranged to visit Carleton---a relatively obscure college from the perspective of a distinguished member of the USSR Academy of Sciences---already contemplated leaving Russia for the stability of the North American continent he had experienced  during his visit to the ICM in 1924.

There is one additional item of evidence one can advance to support this hypothesis.  In the first of his two papers in the \ti{Transactions} of the AMS for 1928 (on the representation of numbers by quadratic forms), referring to his earlier work on this subject in a series of ten papers in Russian journals that had appeared between 1913 and 1926, Uspensky says that he had published these earlier investigations \ti{``so far as it was possible under the circumstances''}.  We take this to be an indirect but clear indication that when Uspensky left for Carleton he had already found the working conditions in revolutionary Russia to have significantly interfered with his ability to do mathematics.

\subsection{Stanford}

Uspensky spent the summer of 1929 teaching in Minnesota,\footnote{Lecturing on number theory in the first term (June 18--July 27) and recent advances in mathematical probability in the second (July 29--August 31), as well supervising reading in advanced mathematics during both;  see  Notes and News, \ti{The American Mathematical Monthly} 36 (3), p.\  421.} and then moved to Stanford that Fall. Initially appointed Acting Professor of Mathematics for the first two years (1929--1931), he became a permanent member of the faculty in 1931.  This presumably came about thanks to the efforts of the chair of the Department of Mathematics,  \ti{Hans Frederik Blichfeldt} (January 9, 1873--November 16, 1945).

Blichfeldt, who was born in Denmark but came to the US in 1888 when his family emigrated, had impeccable mathematical credentials:  after graduating with an AB from Stanford in 1896 he---like many US mathematicians at the time---went to Europe for his doctorate, studying at the University of Leipzig under the great Sophus Lie and receiving his PhD there in 1898.\footnote{The Mathematics Genealogy Project gives the year of Blichfeldt's doctorate as 1900, but this was actually the year his thesis was published in the \ti{American Journal of Mathematics} (Blichfeldt, 1900).}  He then returned to Stanford, where he remained for the rest of his professional life, becoming a Full Professor in 1913, a member of the National Academy of Sciences in 1920, and Department Chair from 1927 to 1938.  He did important work in group theory and number theory (including co-authoring a book with Dickson).  For further information about Blichfeldt, see his obituary in the \ti{Bulletin} of the AMS by Dickson (1947), his National Academy of Sciences biographical memoir by Eric Temple Bell (1951), his entry in the \ti{Dictionary of Scientific Biography} (Miller, 2008), and Roydan's history of the Stanford Mathematics Department (1989, pp. 238--248).

Royden (1989, p.\ 244) relates that after Blichfeldt became chair in 1927, ``the Stanford mathematics department had a steady stream of major mathematicians as visiting faculty, mostly for the summer quarter''.  (The Europeans among these included Harald Bohr, Edmund Landau, and Gabor Szeg\"{o}.)  It was thus not surprising that he seized the opportunity of hiring Uspensky when the latter moved to the US.  This was a major commitment at the time because the Stanford Mathematics Department was then quite small: Royden (1989, p.\ 248)  reports  that in 1938, the year Blichfeldt retired, the Department had only three Full Professors:  Blichfeldt, Uspensky, and W.\ A.\ Manning.\footnote{\ti{William Albert Manning} (1876--1972) was Stanford's first PhD in Mathematics (1904), taught at Stanford for 40 years, and was the father of Laurence Albert Manning (1923--2015), himself later a professor in Stanford's Department of Electrical Engineering for 40 years.}

Up until 1935 Uspensky published a steady stream of research papers, but then abruptly changed the form and direction of his work in several ways.  First, he turned to writing books instead of papers, his books on probability, number theory, and the theory of equations appearing in 1937, 1939, and 1948 (the last finished in late 1946, shortly before his death).  At the same time he developed an interest in applied mathematics:  in the last decade of his life ``he took an active part in Applied Mechanics seminars and presented there many talks of great interest'' (P\'{o}lya, Szeg\"{o}, and Young, 1947, p.\ 2).  His colleague and fellow refugee, the engineer Timoshenko, who would have participated in these, recounted  shortly after Uspensky's death that
\begin{quote}
He participated in a seminar on applied mechanics and gave talks on various topics of applied mathematics, in such a way as to show clearly to engineers the importance of mathematics in treating engineering problems.  The students appreciated these talks very much, and enjoyed also the informal discussions after seminar in which J.\ V.\ usually participated by telling them stories and anecdotes from the biographies and memoirs of famous scientists, about whom he had a great fund of knowledge  (Timoshenko, 1947, p. 6).\footnote{Thus an unpublished remembrance written shortly after Uspensky's death.  Such talks were not always successful however from Timoshenko's perspective:  writing many years later (Timoshenko, 1968, p. 337), he remembered that when he asked Uspensky to lecture on partial differential equations, the lectures were ``purely theoretical'' and not what the engineering students needed.}
\end{quote}
In the last three years of his life Uspensky also taught himself Spanish and published five papers in Spanish during this period. 

Towards the end of 1946 Uspensky suffered from ill health and, although continuing to write and lecture, planned to retire from Stanford at the end of the 1946--47 academic year.  But he suffered a heart attack on Thursday, January 23, 1947 and was taken to the Palo Alto Hospital, where he died the next Monday at 2:30 in the afternoon.\footnote{Some sources give the place of death as San Francisco, but the local obituaries of the time all state he died at the Palo Alto Hospital.  The details of when he was stricken, and  when and where he died, have been drawn from these.}

\section{Uspensky's Students}

Uspensky had students in both Russia and the US.  All of these---although in very different ways---led interesting lives.

\subsection{Students in Russia}

Despite his wide range of mathematical interests, Uspensky was first and foremost a number theorist.  Of the five leading number theorists in Leningrad in the 1920s who did not leave Russia---Delone, Ivanov, Kuzmin, Venkov, and Vinogradov (Demidov, 2015, p.\ 89)---three were students of Uspensky.\footnote{The course of graduate study in Russia during the 1918--1934 period  is complicated.  On October 1, 1918, the Council of People's Commissars abolished all academic ranks and advanced degrees.  Instead, advanced training was initially supported by ``professorial scholarships'' for students  preparing for an academic career.  After 1925 the process of selection and mentoring for such studies was formalized and participants termed ``aspirants'' (\ti{aspirantov}).   By a decree of September 19, 1932, the standards for such programs were tightened and in research institutes the writing of a dissertation became required after a two to three year program of study.  Finally, on January 13, 1934 academic ranks and advanced degrees were reinstated. See DeWitt (1961, p.\ 422), on which the information in this paragraph is based.}\\

\ti{Ivan Matveevich Vinogradov} (September 14, 1891--March 20, 1983) was both Uspensky's first and best-known student.
Vinogradov was a central figure in modern analytic number theory, a member of the USSR Academy of Sciences, the Director of the Steklov Institute of Mathematics in Moscow for half a century (from 1934 until his death), and a Fellow of the Royal Society of London.   His best known result is that every sufficiently large odd number is a sum of three primes.  His \ti{Foundations of the Theory of Numbers}, first published in 1936, is still in print today in an English translation of the 5th edition.   He is too well-known to require discussion here;  see his Royal Society biographical memoir (Cassels and Vaughn, 1985) and Karatsuba (1981).\footnote{The Mathematical Genealogy Project entry for Vinogradov, while listing Uspensky as his advisor, does not list a degree, year when it was awarded, or title of the dissertation.  Presumably this was because of the abolition of degrees noted in the preceding footnote.}\\

We have already encountered \ti{Rodion Osievich Kuzmin} (October 9, 1891--March 24, 1949):  his 1928 solution of Gauss's challenge problem to Laplace was given in the third  appendix of Uspensky's  \ti{Introduction to Mathematical Probability}.  After completing his undergraduate studies at the University of Petrograd in 1916, Kuzmin continued on to do  graduate work there, but was forced to leave in 1918 due to wartime conditions.  It was during this short period of time that he would have studied under Uspensky.  Like Vinogradov, Kuzmin moved to the University of Perm, where he stayed until 1922, returning then to Petrograd and its Polytechnic Institute, where he remained for the rest of his life.  In 1935 he received his \ti{doktorat}, shortly after this degree was reinstated by the Soviet government.\footnote{The Mathematics Genealogy Project states Kuzmin received the degree of Doctor of Sciences in 1935 and identifies his advisor as Uspensky, but the actual circumstances seem unclear:  Kuzmin had already been teaching for many years in 1935 and Uspensky had left Russia six years earlier.  The 1935 degree may have been awarded on the basis of a dissertation drawing on prior work but only submitted after the reinstatement of advanced degrees, but even so it seems highly improbable that at the onset of the Great Terror an \'{e}migr\'{e} who had resigned from the USSR Academy of Sciences would have been listed as Kuzmin's advisor.}
In addition to Gauss's problem, Kuzmin proved in 1930 that if $a$ is algebraic and $b$ is a real quadratic irrational, then $a^b$ is transcendental.  This answered a then open problem, the nature of the  number $2^{\sqrt{2}}$ (the so-called Gelfond-Schneider constant).  This was a special case of the celebrated Gelfond-Schneider theorem (proved in 1937), which assumes only that $b$ is irrational.  Other notable work include his 1927 paper on the \ti{Kuzmin-Landau inequality}.  For further information see Venkov and Natanson (1949), Sviderskaya (2002).\\

\ti{Boris Alekseevich Venkov} (July 31, 1900--December 13, 1962) worked with Uspensky while still a student;  their work was later reported on in a joint paper presented to the 1924 International Congress of Mathematicians in 1924.  Venkov graduated from the Leningrad State University in 1925 and became a professor there in 1935;  his work was in the field of the theory of numbers. His contributions include  an elementary proof of the Dirichlet formulas for a number of the classes of binary squared forms (1928), a theory of reduction for positive-definite quadratic forms (1940), and a characterization of polyhedra (1954 and 1959).  \ti{Venkov graphs} are named after him.  For further information, see Malyshev and Faddeev (1961).

\subsection{Students at Stanford}

Uspensky supervised the doctoral dissertations of  five students at Stanford.  The first of these, \ti{Harold Maile Bacon} (January 13, 1907--August 22, 1992), received his PhD from Stanford in 1933;  his thesis (``An Extension of Kronecker's Theorem'') on a topic in number theory.  He spent his entire career at Stanford, primarily as a long-time instructor of calculus courses at Stanford, and was known to generations of Stanford undergraduates.

Bacon was a ``hand-off'' to Uspensky from Harald Bohr.  Years later Bacon related the following charming story as to how this came to be:

\begin{quote}
Bohr was a very kind man.  I remember my being in Professor Blichfeldt's office shortly after I returned to Stanford in 1930 to continue my graduate work after my master's degree and a year's absence working for an insurance company under the mistaken impression that I wanted to become an actuary. Blichfeldt and I were discussing my getting started on work that might lead to a dissertation. Just then Bohr came into the office. Blichfeldt turned to him and, indicating me, said `Here's a man who is looking for a thesis topic.  How would you like to suggest one, and be his adviser?'  Bohr bowed, smiled and very courteously replied, `I should be honored.'  He generously acted as my supervisor for the remainder of the year he was at Stanford. When he left, I was most fortunate to have Uspensky take over and see me through to the completion of my work on the dissertation. It was indeed a great privilege to have two such inspiring men as my friends and advisers at the beginning of my career.	[Royden (1989, pp.\ 245--6).]
\end{quote}

Because he was such a visible figure at Stanford, a considerable amount has been written about Bacon;  see Royden (1989, pp.\ 244--247), Royden et al.\ (1992), Stanford (1992), Jellison (1997), Albers and Alexanderson (2011, Chapter 3).

\bigskip

\ti{Maxwell Alfred Heaslet} (February 17, 1907--July 13, 1976), Uspensky's next student at Stanford, received his PhD there in 1934 (dissertation ``Concerning the Development Coefficients of an Aequianharmonic Function'').  Information about him is surprisingly scanty.  Here is what we have been able to find, piecing together several scattered sources.

Born in Bentonville, Arkansas, the son of Walter Monroe Heaslet and Nancy Angeline Austin, Heaslet received both his bachelor's and master's degrees from the University of Oklahoma before moving to Stanford in 1928.  After marrying Helen Virginia Camp in June 1935, he joined the faculty of San Jose State College (now University), where he taught mathematics and physics for seven years.  He began working at the Ames Research Center (housed in the Moffett Federal Airfield in nearby Mountain View) in 1942.  He appears to have been initially on loan from San Jose, because he is listed in a departmental history (Jackson, n.\ d.) as being a member of the faculty from 1935 to 1945.  If so, this arrangement presumably ceased at the end of the war, and Heaslet elected to stay on at Ames;  he eventually became head of the theoretical aerodynamics branch of the Theoretical and Applied Research Division.  He often collaborated with Harvard Lomax (1922--1999), who had a distinguished career in aeronautics and computational fluid dynamics (Seabass, 2002), and wrote many papers during his time at Ames;  see Hartmann (1970).
 Due to ill health Heaslet was forced to retire from Ames in 1959 and moved to Florida where he remained for the rest of his life.  He had two children, Austin and Jonathan Heaslet, and a brother, Walter.\\

From a purely academic standpoint Uspensky's most successful student  was \ti{Orville Goodwin Harrold, Jr.} (September 2, 1909--May 16, 1988);  Harrold, who received his PhD in 1936 (dissertation ``On the Expansion of the Remainder in the Open Type Newton-Cotes Quadratic Formula''), had 14 students and 135 descendants.  After the war he worked in topology, teaching at the University of Tennessee at Knoxville for more than a decade, and after that Florida State University, where he was the chair from 1964 to 1974;  the Orville G. Harrold Professor of Mathematics is named in his honor.  He was elected a Fellow of the AAAS in 1954, was a Guggenheim Fellow in 1957-58, member of the Institute for Advanced Studies twice (1958 and 1964), and an Associate Secretary of the AMS from 1965 to 1976.  He was married to Gladys E.\ Buell on June 30, 1934, and had a son Jeffrey Buell Harrold. \\

At the other extreme, \ti{Franklin Alfred Butter, Jr.\ }(February 1, 1910--November 27, 1972), who received all three of his degrees from Stanford (AB 1930, AM 1931, PhD 1936) and graduated the same year as Harrold (the title of his dissertation ``A Contribution to the Theory of the Arithmetic-Geometric Mean''), is more elusive.    We know the basic details of his career up to 1946  thanks to his entry in the Stanford University \ti{Annual Register} for 1945--46, when he temporarily returned to Stanford as an Acting Assistant Professor.\footnote{\ti{Stanford University Fifty-fifth Annual Register 1945--46.  Stanford University, Published by the University}.   An invaluable historical resource, the registers are available online at \url{https://exhibits.stanford.edu/stanford-pubs/browse/annual-register-1891-1947} (last accessed January 30, 2022).  
}  
After he received his PhD at Stanford, he spent the year 1936--1937 at Washington University in St.\ Louis as a Rockefeller Science Fund Research Fellow, collaborating with Szeg\"{o} in the preparation of Szeg\"{o}'s classic book \ti{Orthogonal Polynomials} (Szeg\"{o}, 1938, p.\ vii), and then went on to USC for six years.  He appears to have published nothing during his time at USC, however, leaving there in 1943 and spending the next three years in a series of one-year positions at the University of Wisconsin, Lawrence College, and Stanford.  But then his life took an unexpected turn: in 1942--43 he had been a mathematical consultant in the Engineering Department at the Douglas Aircraft Company, and this appears to have generated an interest in aerospace science (see Butter, 1945, a book review of \ti{Aircraft Analytic Geometry}).  In any case, he left academia in 1946 to work in the aerospace industry:  first as a mathematician and later Research Physicist at the Hughes Aircraft Company, Culver City, California, 1946--57; then an engineering specialist at Northrop Aircraft, Inc., Hawthorne, California, 1957--1961;  and finally a Staff Engineer at the Aerospace Corporation, El Segundo, California, 1961--65 (all three of these located in Los Angeles County).  But he returned to teaching in 1965, when he was appointed an Associate Professor at California State College at Long Beach, and promoted to Professor there in 1970.\footnote{The details of Butter's life from 1946 on have been teased from the \ti{News and Notices} pages of \ti{The American Mathematical Monthly},  
53 (1946), p.\ 603;  65 (1958), p.\ 63;  73 (1966), p.\ 108;  78 (1971), p.\ 107;  and an announcement in the \ti{Notices of the AMS},  8 (1961), p.\ 486.}\\

Finally there is \ti{Carl Douglas Olds} (May 11, 1912--November 11, 1979), \, Uspensky's last student, who received his PhD in 1943, his dissertation, like Bacon's, in number theory (``On the Number of Representations of the Square of an Integer as the Sum of an Odd Number of Squares'').  He was an acting instructor at Stanford from 1935 to 1940, and an assistant professor at Purdue from 1940 to 1945, before joining the faculty of San Jose State in 1945 (perhaps as a replacement for Heaslet), where he remained for the rest of his career.  He was an active member of the Mathematical Association of America, winning the 1973 Chauvenet prize for his article on the continued fraction expansion of $e$ (Olds, 1970), and may be familiar to some because of his lovely book \ti{Continued Fractions} (1963).

\section{Bibliography of J.\ V.\ Uspensky}

AJM:  \ti{American Journal of Mathematics}\\

AMM:  \ti{The American Mathematical Monthly}\\

BAMS:  \ti{Bulletin of the American Mathematical Society}\\

TAMS:  \ti{Transactions of the American Mathematical Society}
\bigskip

\subsection{Thesis}

\nin 1910  Some applications of continuous parameters to number theory.  Online HathiTrust copy at https://catalog.hathitrust.org/Record/006087199  (accessed March 23, 2020).\\

Studies minimization problems for pairs of linear forms $\alpha x + \beta y$ and $\gamma x + \delta y$, continuing earlier research of Voronoi and Minkowski. The main result of the thesis is an algorithm for the reduction of forms of the type $\alpha x^2 + 2 b x y + c y^2$ that depend on $i$ and $e^{2\pi i/3}$.  [See Markov et al.\ (1921, p.\ 4).]\\

\subsection{Books}

\nin 1937  \ti{Introduction to Mathematical Probability}.  New York:  McGraw-Hill.\\

\nin \tb{Reviews}:  F.\ N.\ David, \ti{Biometrika} 30, 194--195;  J\ A.\ Greenwood,
\ti{The American Mathematical Monthly} 45, 471;  H.\ T.\ H.\ Piaggio, \ti{The Mathematical Gazette} 22, 202--204;  Anonymous, \ti{Nature} 141, 769.\\

\bigskip

\nin 1939  \ti{Elementary Number Theory} (co-author Maxwell\ A.\ Heaslet).  New York: McGraw-Hill. \\

\nin \tb{Reviews}:   R.\ Oldenburger, \ti{Bull.\ Amer.\ Math.\ Soc.} 46, 202--205;  G.\ F.\ Cramer, \ti{National Mathematics Magazine} 14, 494;  L.\ J.\ Mordell, \ti{The Mathematical Gazette} 24, 295--298;  H.\ Davenport, \ti{Nature} 146, 418--419.  MR0000236  (Reviewer: D. H. Lehmer).\\

\bigskip

\nin 1948 \ti{Theory of Equations}.  New York:  McGraw-Hill.\\

\nin \tb{Reviews}:  Garrett Birkhoff, \ti{Science} 109, 429;  T.\ A.\ Brown, \ti{The Mathematical Gazette} 34, 140--142;  L.\ E.\ Bush, \ti{The American Mathematical Monthly} 56, 348--350;  Kenneth May, \ti{Popular Astronomy} 57, 46--47.\\

\subsection{Papers}

Unless otherwise noted, all papers up to 1925 are in Russian and all papers from 1927 on are in English.\\

\nin 1906 \\

\nin 1.\ On whole numbers formed with the 5th root of unity.  \ti{Matematicheskii Sbornik} 26, 1--17.  [Proves $Z[\zeta_5]$ is Euclidean.]\\

\nin 1909 \\

\nin 1.\ Note sur les nombres entiers d\'{e}pendant d'une racine cinqui\`{e}me de l'unit\'{e}. \ti{Math. Ann.} 66, 109--112 (in French).  \\

\nin 1912 \\

\nin 1.  Sur une s\'{e}rie asymptotique d'Euler.  \ti{Archiv der Mathematik und Physik} 19, 370--1 (in French).\\
 
\nin 1913  \\

\nin 1.\  Arithmetical proof of the Kronecker relations between the class numbers of binary quadratic forms.  \ti{Matematicheskii Sbornik} 29, 26--52.  [Provides an elementary proof of the Kronecker relations.  This line of investigation was later continued by Uspensky in his seven papers in the years 1925--1926.]\\

\nin 2.\  On arithmetical theorems given by Stieltjes. \ti{Bulletin of the Mathematical Society in Kharkov} 14, 7--30.\\

\nin 3.\   On the representation of numbers by sums of squares. \ti{Bulletin of the Mathematical Society in Kharkov} 14, 31--64.\\

\nin 4.\  On certain arithmetic theorems, I.  \ti{Bulletin of the Mathematical Society in Kharkov} 14, 88--96.\\

\nin 1914 \\

\nin 1.\ On certain arithmetic theorems, II.  \ti{Bulletin of the Mathematical Society in Kharkov} 14. 97--99.\\

\nin 2.\ On the possibility of representing prime numbers by some of the simplest quadratic forms.  Kazan.  [Cited in Markov et al., 1921, p.\ 5;  no further bibliographic details given.]\\

\nin 3.\  A rule for determining the sign in the equality $1 \cdot 2 \cdot \dots \frac{p-1}{2} \equiv \pm 1 \pmod{p}$ for a prime $p$ of the form $4 \mu + 3$.  Kazan.  [Cited in Markov et al., 1921, p.\ 6;  no further bibliographic details given.]\\

\nin 1915  \\

\nin 1.\ On the class numbers of positive Hermite forms.  \ti{Bulletin of the Imperial Academy of Sciences of St.-Petersburg} 9, 1769--1800.\\
\newline

\nin 1916  \\

\nin 1.\ On the representation of numbers by the quadratic forms with 4 and 6 variables, I.  \ti{Bulletin of the Mathematical Society in Kharkov} 5, 81--112.\\

\nin 2.\  On the representation of numbers by the quadratic forms with 4 and 6 variables, II.  \ti{Bulletin of the Mathematical Society in Kharkov} 5, 113--147.\\

\nin 3.\   On the convergence of quadrature formulas between infinite limits.  \ti{Bulletin of the Imperial Academy of Sciences of St.-Petersburg} 10, 851--866.  [Extended in Uspensky, 1928/2]\\

\nin 4.\   On the development of functions in series arising from the polynomials $e^x \frac{d^n x^n e^{-x}}{dx^n}$.  \ti{Bulletin of the Imperial Academy of Sciences of St.-Petersburg} 10, 1173--1202.\\

\nin 1920  \\

\nin 1.\ Asymptotic formulae for numerical functions which occur in the theory of the partition of numbers into summands.  \ti{Bulletin of the Academy of Sciences of Russia} 14,  199--218. \\

\nin 1921  \\

\nin 1.\ About approximate expressions for the coefficients of distant terms in the development of the equation of the center into a series by the sine of multiples of the mean anomaly.  \ti{Bulletin of the Academy of Sciences of Russia} 15,  333--342.\\

\nin 1923  \\

\nin 1.\ Note on the scientific work of A.\ A.\ Markov.  \ti{Bulletin of the Academy of Sciences of Russia} 17, 19--34.\\

\nin 1924  \\

\nin 1.\ On a problem of Jean Bernoulli.  \ti{Bulletin of the Academy of Sciences of Russia} 18, 67--84.\\

\nin 1925  \\

\nin  1.\ Sur les valeurs asymptotiques des coefficients de Cotes (in French).  BAMS 31, 145--156.  [Describes the asymptotic behavior of the coefficients in the Newton-Cotes method.  ``Attention should be called to the fact that many formulas on p. 147 \dots are marred by typographical errors'' (Uspensky, 1935/1, p.\ 382).]\\

\nin 2.\ Note sur le nombre des repr\'{e}sentations des nombres par une somme d'un nombre pair de carr\'{e}s.  \ti{Bulletin of the Academy of Sciences of the USSR} 19, 647--642.\\

\nin 3.\ Sur les relations entre les nombres des classes des formes quadratiques binaires et positives.  Premier M\'{e}moire, I.  \ti{Bulletin of the Academy of Sciences of the USSR} 19, 599--620.\\

\nin 4.\ Sur les relations entre les nombres des classes des formes quadratiques binaires et positives.  Premier M\'{e}moire, II.  \ti{Bulletin of the Academy of Sciences of the USSR} 19, 763--784.\\

\nin 1926  \\

\nin 1.\ Sur les relations entre les nombres des classes des formes quadratiques binaires et positives.  Deuxi\`{e}me M\'{e}moire, I.  \ti{Bulletin of the Academy of Sciences of the USSR} 20, 25--38.\\

\nin 2.\ Sur les relations entre les nombres des classes des formes quadratiques binaires et positives.  Deuxi\`{e}me M\'{e}moire, II.  \ti{Bulletin of the Academy of Sciences of the USSR} 20, 175--196.\\

\nin 3.\ Sur les relations entre les nombres des classes des formes quadratiques binaires et positives.  Troisi\`{e}me M\'{e}moire.  \ti{Bulletin of the Academy of Sciences of the USSR} 20, 327--348.\\

\nin 4.\ Sur les relations entre les nombres des classes des formes quadratiques binaires et positives.  Quatri\`{e}me M\'{e}moire.  \ti{Bulletin of the Academy of Sciences of the USSR} 20, 547--566.\\

\nin 5.\ Sur les relations entre les nombres des classes des formes quadratiques binaires et positives.  Cinqui\`{e}me M\'{e}moire.  \ti{Bulletin of the Academy of Sciences of the USSR}, 20 619--642.\\

The last seven papers rigorously derive class-number relations using elementary methods;  these relations were already known but had been obtained previously using the more complex machinery of elliptic functions.  For later discussion, see Dwyer (1937)\\

\nin 6.\ Vladimir Andreevich Steklov.  \ti{Bulletin of the Academy of Sciences of the USSR} 20, 10--11, 837--856.\\

\nin 1927  \\

\nin 1.\ Note on the computation of roots.  AMM 34, 130--134.\\

\nin 2.\ A curious case of the use of mathematical induction in geometry.  AMM  34, 247--250.)\\

\nin 3.\ On a problem arising out of the theory of a certain game.  AMM 34, 516--521.\\

\nin 4.\ On the development of arbitrary functions in series of Hermite's and Laguerre's polynomials.  \ti{Annals of Mathematics} 28, 593--619. [Uspensky says in the introduction that his paper 1916/4 ``was written in Russian and published during the time when all the international relations were broken, and this may sufficiently account for the fact that it could pass unnoticed. However, as the method used by the author can successfully be applied whenever asymptotic expressions of a certain type exist, it seems worth while to  reproduce the essential parts of this paper in a modified and simplified form''.  See Sansone (1950).]  \\

\nin 1928  \\

\nin 1.\ On Jacobi's arithmetical theorems concerning the simultaneous representation of numbers by two different quadratic forms.  TAMS 30, 385--404.  [``Uspensky developed the elementary methods which seem to have been used by Liouville in a series of papers  published in Russian;  references will be found in [this Transactions paper].  He carries his analysis up to $2s = 12$, and states that his methods enable him to prove Boulyguine's general formulae (Hardy and Wright,1960, p.\ 316).]\\

\nin 2.\ On the convergence of quadrature formulas related to an infinite interval.  TAMS 30, 542--559.  [$L_2$ convergence of Lagrange interpolation on the entire real line.  Nevai (1986, p.\ 118) notes Uspensky's priority in dealing with the problem more than three decades earlier than other literature.]   \\

\nin 3.\ On Gierster's classnumber relations.  AJM 50, 93--122.  [Applies results in his 1925-26 series of memoirs.]\\

\nin 4.\ On some new class-number relations.  \ti{Proceedings of the International Mathematical Congress} (J. C. Fields, ed.), Vol. 1, 315-317.  [Written with his student Boris Venkov]\\
\newline

\nin 1929  \\

\nin 1.\ On the number of representations of integers by certain ternary quadratic forms.  AJM  51, 51--60.  [See Kaplansky, 2013, for later context.] \\

\nin 1930  \\

\nin 1.\ On the reduction of the indefinite binary quadratic forms.  BAMS 36, 710--718.\\

\nin 2.\  On incomplete numerical functions.  BAMS 36(10), 743--754.  [Applies results in his 1925-26 series of memoirs.]\\

\nin 1931  \\

\nin 1.\ A method for finding units in cubic orders of a negative discriminant.  TAMS 33, 1--22.  [An extension of  the results in Zolotarev's classic paper of 1869.] \\

\nin 2.\ On Ch. Jordan's series for probability.  \ti{Annals of Mathematics} 32(2), 306--312. [Uspensky's paper ``has been overlooked by most authors dealing with [Poisson] approximation problems'' (Deheuvels and Pfeifer, 1988, pp.\ 671--672);  Deheuvels and Pfeiffer go on to establish a relationship between Uspensky's approximation and the modern Poisson convolution semigroup approach, and exploit this to simplify and extend prior results in the literature.] \\

\nin 1932  \\

\nin 1.\ On the problem of runs.  AMM 39, 322--326.\\

\nin 1933 \\

\nin 1.\ A minimum problem.  BAMS 40, 5--10.\\

\nin 1934  \\

\nin 1.\ On an expansion of the remainder in the Gaussian quadrature formula.  BAMS 49, 871--876.  [``In this article I shall prove that the remainder in the Gauss\-ian formula can be expanded into a series possessing all the properties of the classical Euler-Maclaurin expansion.'']\\

\nin 1935  \\

\nin 1.\ On the expansion of the remainder in the Newton-Cotes formula.  TAMS 37, 381--396.  [Extends work in 1934/1]\\

\nin 1944  \\

\nin 1.\ Elementary derivation of the series for $\sin x$ and $\cos x$.  (Spanish)  \ti{Math. Notae} 4, 1--10.\\

\nin 2.\   A new proof of Jacobi's theorem.  (Spanish)  \ti{Math. Notae} 4, 80--89.\\

\nin 1945  \\

\nin 1.\ On the problem of the ruin of gamblers.  (Spanish)  \ti{Publ. Inst. Mat. Univ. Nac. Litoral } 7, 155--186.\\

\nin 2.\ Sur la m\'{e}thode de Laplace dans la th\'{e}orie de l'attraction des ellipsoides homog\`{e}nes. (French) \ti{Publ. Inst. Mat. Univ. Nac. Litoral} 5, 63--71.\\

\nin 3.\ On the arithmetico-geometric means of Gauss, I--III.  (Spanish)  \ti{Math. Notae} 5, 1--28, 57--88, 129--161.\\

\nin 1946  \\

\nin 1.\ On a problem of John Bernoulli, I--IV. (Spanish) \ti{Rev. Un. Mat. Argentina} 11, 141--154, 165--183, 293--255;  12, 10--19.\\

\subsection{Book Reviews and Letters}

\nin 1932  Review of \ti{Studies in the Theory of Numbers}, by L.\ E.\ Dickson. BAMS 38, 463--465.\\

\nin 1940  Review of \ti{Development of the Minkowski Geometry of Numbers} by Harris Hancock. \ti{National Mathematics Magazine} 14, 423--424.\\

\nin 1944 Remarks on the History of Science in Russia.  \ti{Science} 100, 193--194 (with S.\ P.\ Timoshenko).\\

\nin 1946  Book Review:  P. L. Chebyshev, \ti{Collected works}. Vol.\ 1. \ti{Theory of numbers}, BAMS 52, 50.\\

\subsection{Problems and Solutions}

Uspensky was a frequent contributor to the Problem Section of the \ti{American Mathematical Monthly}, submitting at least eleven problems and three solutions to it:  Number 3251, 34:4 (Apr.\ 1927), p.\ 216 [solution given by Uspensky 37:6 (Jun.--Jul.\ 1930), pp.\ 318--319];  Number 3278,  34:7 (Aug.--Sep.\ 1927), p.\ 381;   Number 3290, 34:9 (Nov.\ 1927), p.\ 491 [reproposed by C.\ D.\ Olds -- Uspensky's student -- as Number 4400, 57:6 (June 1950), p.\ 420];    Number 3304, 35:1 (Jan.\ 1928), p.\ 41; Number 3312, 35:3 (Mar. 1928), p.\ 154;  Number 3343, 35:8 (Oct.\ 1928), p.\ 446 [solved by H.\ Langman, 36:8 (Oct.\ 1929), p.\ 450];    Number 3354, 35:10 (Dec.\ 1928), p.\ 563;   Number 3389, 36:8 (Oct.\ 1929), p.\ 448;  Number 3408, 37:1 (Jan.\ 1930), p.\ 38 [solution given by Uspensky 39:3 (Mar.\ 1932), pp.\ 176--177];  Number 3534, 39:2 (Feb.\ 1932), p.\ 116 [but see 68:8 (Oct.\ 1961), p.\ 814 indicating an error in the formula proposed by Uspensky, and 69:2 (Feb.\ 1962), pp.\ 172--173, showing that this was due to a simple typographical error];  solution to Number 3588 in 40:10 (Dec.\ 1933), p.\ 614 [note ``Solved also by L. Zander -- his wife!]\\

\subsection{Other writings}  In addition to we list above, the 1921 report by Markov et al.\  to the Russian Academy of Sciences lists two papers whose publication status was unknown to them:\\

\nin a.  Application of the Poisson summation method to some types of Chebyshev polynomials.\\

\nin [Markov et al.\ (1921, p.\ 6) say:  ``This paper should have been published in the jubilee edition of the Notices (Izvestiya) of the Kazan Mathematical Society. It is unknown what happened to it.''  It may be this is the paper Uspensky (1927/5, p.\ 593) is referring to when, commenting on Hille (1926), he wrote: ``From the same point of view, but in a slightly different way, the subject was treated in a paper by the author presented to the Math.\ Society in Kazan in 1916. Owing, however, to the political disturbances of that time, this paper never was published and was finally lost''.]\\

\nin b.  On the approximate expression of the deleted terms in the expansion of the equation center in a row of sines multiple of the average anomaly.  Markov et al.\ say:  ``This paper should have been published in the ``Notices'' (\ti{Izvestiya}) of the Perm State University, but it is unknown whether it has been printed by now.''\\

At the end of their report Markov et al.\ add ``Besides  the aforementioned works there are preprints for the following'', and then list without date what appear to be five pedagogical books or pamphlets written by Uspensky:\\

1.  Lectures on non-Euclidean geometry. 126 p.\\

2.  Theory of biquadratic and cubic residues. 160 p.\\

3.  Continuous fractions and their applications. 160 p.\\

4.   Introduction to additive number theory. 80 p.\\

5.  Liouville methods in number theory. 250 p.\\

\paragraph{Translations}

\nin J.\ Bernoulli, \ti{Ars Conjectandi, Part 4}, translated from the French by Ya. V. Upenskii, and with an introduction by A.\ A.\ Markov.
St.-Petersburg: Imp. Akademiya Nauk, 1913.  Reprinted with additional commentary by Yu.\ V.\ Prokhorov, 1986.\\

\paragraph{Books in Russian}

\nin \ti{Ocherk istorii logarifmov} (Essay on the History of Logarithms).  Petrograd:  Nauchnoe knigoizdatel'stvo, 1923.\\

\section{Bibliography}

\nin Albers, Donald J. and Alexanderson, Gerald L., eds. (2011).  \ti{Fascinating Mathematical People:  Interviews and Memoirs}, Chapter 3, ``Harold M. Bacon''.  Princeton University Press.\\

\nin Bahadur, R.\ R,\ (1960).  Some approximations to the binomial distribution function.  \ti{Annals of Mathematical Statistics} 31, 43--54.\\

\nin Bahadur, R.\ R.\ (1966).  A note on quantiles in large samples.  \ti{Annals of Mathematical Statistics} 37, 577--580.\\

\nin Beatty, S. (1926).  Problem 3173.  \ti{American Mathematical Monthly} 33, 159.\\

\nin Beatty, S.\ et al. (1927).  Solutions to Problem 3173.  \ti{American Mathematical Monthly} 34, 159--160.\\

\nin Bell, Eric Temple (1951).  Hans Frederick Blichfeldt 1873--1945.  \ti{Biographical Memoirs, V.\ 26}, National Academy of Sciences, Washington, DC., 181--189.\\

\nin Bennett, George (1962).  Probability inequalities for the sum of independent random variables.  \ti{Journal of the American Statistical Association} 57, 33--45.\\

\nin Bernoulli, Jacob (2006).  \ti{The Art of Conjecturing} (translated, with an introduction and notes by Edith Dudley Sylla).  Baltimore:  Johns Hopkins University Press.\\

\nin Bernstein, Sergei N.\ (1924) . On a variant of Chebyshev's inequality and on the accuracy of Laplace's formula (in Russian). \ti{Uchenie Zapiski Nauchno-Issledovatel'nikh Kafedr Ukrainy, Otdelenie Matematiki} 1, 38--48.\\

\nin Bernstein, Sergei N.\ (1927).  \ti{Teoriya Veroiatnostei}  [\ti{The Theory of Probability}]. Moscow--Leningrad: Gosudarstvennoe Izdatel'stvo.\\  

\nin Bertrand, J.\ (1888).  \ti{Calcul des probabilit\'{e}s}.  Paris:  Gauthier-Villars (reprinted 1889, 1907).  Reprinted, 1972;  New York:  Chelsea.\\

\nin Bhatia, R.\ A.\ (2008).   Conversation with S.\ R.\ S.\ Varadhan. \ti{The Mathematical Intelligencer} 30, 24--42.\\

\nin Blackwell, David (1954).  On optimal systems.  \ti{The Annals of Mathematical Statistics} 25, 394--397.\\

\nin Blichfeldt, Hans Frederick (1900).  On a certain class of groups of transformation in space of three dimensions.  \ti{American Journal of Mathematics}  22, 113--120.\\

\nin Borwein, J.\ M.\, Choi, S.\ K.\, and Pigulla, W.\ (2005).  Continued fractions of tails of hypergeometric series.  \ti{American Mathematical Monthly} 112, 494--501.\\

\nin Boucheron, S., Lugosi, G., and Massart, P.\ (2012).  \ti{Concentration Inequalities:  A Nonasymptotic Theory of Independence}.  Oxford:  Clarendon Press.\\

\nin Butter, F.\ A.\ (1945).  Book review of \ti{Aircraft Analytic Geometry}, by J.\ J.\ Apalategui and L.\ J.\ Adams.  \ti{The American Mathematical Monthly} 52, 338--340.\\

\nin Cassels, J.\ W.\ S.\ and Vaughan, R.\ C.\ (1992).   Ivan Matveevich Vinogradov. \ti{Biographical Memoirs of Fellows of the Royal Society} 31, 612--631;  reprinted in the \ti{Bulletin of the London Mathematics Society} 17, p584--600.\\

\nin Chamberlain, Lesley (2006).  \ti{Lenin's Private War}.  New York:  St.\ Martin's Press.\\

\nin Cochran, William G.\ (1952).  The $\chi^2$ test of goodness of fit.  \ti{Annals of Mathematical Statistics} 23, 315--345.\\

\nin Coolidge, Julian Lowell (1925).  \ti{An Introduction to Mathematical Probability}.  Oxford:  Clarendon Press.\\

\nin Craig, Cecil C.\ (1933).  On the Tchebychef inequality of Bernstein.  \ti{Annals of Mathematical Statistics}
 4, 94--102.\\
 
 \nin Cram\'{e}r, H.\ (1976). Half of a century of probability theory: some personal recollections. \ti{Annals of
Probability} 4, 509-546.\\

\nin Deheuvels, P.\ and Pfeifer, D.\ (1988).  On a relationship between Uspensky's theorem and Poisson approximations.  \ti{Annals of the Institute of Statistical Mathematics} 40, 671--681.\\

\nin Demidov, S.\ S. (2015). World War I and mathematics in ``the Russian world''. \ti{Czasopismo Techniczne} 112, 77--92.\\

\nin DeSalvo, Stephen (2021).  Will the real Hardy-Ramanujan formula please stand up?  \ti{Integers} 21, \#A88.\\

\nin DeWitt, Nicholas (1961).  \ti{Education and Professional Employment in the U.S.S.R.}  Washington, D.\ C.:  National Science Foundation. \\

\nin Diaconis, P.\ , Graham, R.\ L.\, and Kantor, W.\  M.\ (1983).  The mathematics of perfect shuffles.  \ti{Advances in Applied Mathematics} 4, 175--196. \\

\nin Diaconis, P.\ and Graham, R.\ L.\ (2012).  \ti{Magical Mathematics:  The Mathematical Ideas That Animate Great Magic Tricks}.  Princeton University Press. \\

\nin Diaconis, P.\ and Zabell, S.\ L.\ (1991).  Closed form summation for classical distributions:  variations on a theme of De Moivre.  \ti{Statistical Science} 6, 284--302.\\

\nin Dickson, L. E. (1947), Hans Frederik Blichfeldt:  1873--1945. \ti{Bulletin of the American Mathematical Society} 53, 882-883.\\

\nin Dudley, R.\ M.\ (1987).  Some inequalities for continued fractions.  \ti{Mathematics of Computation} 49, 585--593.\\

\nin Dutka, Jacques (1981).  The incomplete beta function---a historical profile.  \ti{Archive for History of Exact Sciences} 24, 11--29.\\

\nin  Dwyer, W.\ A.\ (1937).  On certain fundamental identities due to Uspensky.  \ti{American Journal of Mathematics} 59, 290--294.\\

\nin Ermolaeva, N.\ (1997).  Uspensky Yakov Viktorovich.  In \ti{Russian Diaspora: The Golden Book of Emigration, the First Third of the 20th Century, An Encyclopedic Biographical Dictionary}  (edited by V.\ V.\ Shelokhaeva;  executive editor N.\ I.\  Kanishcheva).  Moscow: ROSSPEN.  Available online in Russian at \url{http://www.apmath.spbu.ru/ru/misc/uspenskii.html} (accessed January 5, 2021).\\

\nin Ethier, Stewart  N. (2010).  \ti{The Doctrine of Chances}.  New York:  Springer-Verlag.\\

\nin Ethier, Stewart and Khoshnivisan, Davar (2002).  Bounds on gambler's ruin probabilities in terms of moments.  \ti{Methodology and Computing in Applied Probability} 4, 55--68.\\

\nin Feller, William (1950).  \ti{An Introduction to Probability Theory and its Applications}, Volume 1.  New York:  Wiley.  2nd ed. 1957, 3rd ed. 1968.\\

\nin Finkel, Stuart (2003).  Purging the public intellectual:   the 1922 expulsions from Soviet Russia.  \ti{The Russia Review} 62, 589--613.\\

\nin Fisher, Arne (1922).  \ti{The Mathematical Theory of Probabilities \& Its Application to Frequency Curves \& Statistical Methods}, 2nd ed.  New York:  Macmillan.\\

\nin Fisher, R.\ A.\ (1925).  \ti{Statistical Methods for Research Workers}, 1st ed.  Edinburgh:  Oliver \& Boyd.  14th (posthumous) ed., 1970, New York:  Hafner.\\

\nin Forsyth, C.\ H.\ (1924).  Simple derivations of the formulas for the dispersion of a statistical series.  \ti{The American Mathematical Monthly} 31, 190--196.\\

\nin Gilliland, D., Levental, S., and Xiao, Y.\ (2007).  A note on absorption probabilities in one-dimensional random walk via complex-valued martingales.  \ti{Statistics and Probability Letters} 77, 1098--1105.\\

\nin Graham, R.\ L.\ (1963).  On a theorem of Uspensky.  \ti{The American Mathematical Monthly} 70, 407--409.  [Discusses AMM 1927 516 paper on a certain game]\\

\nin Graham, R.\ L., Lin, S., and Lin, C.-S.\ (1978).  Spectra of numbers.  \ti{Mathematics Magazine} 51, 174--176.\\

\nin Graham, Ron and O'Bryant, Kevin (2005).  A discrete Fourier kernel and Fraenkel's tiling conjecture.  \ti{Acta Arithmetica} 118, 283--304.\\

\nin Hald, Anders (1990).  \ti{A History of Probability and Statistics and Their Applications before 1750}.  New York:  John Wiley \& Sons.\\

\nin Hald, Anders (1998).  \ti{A History of Mathematical Statistics from 1750 to 1930}.  New York:  John Wiley \& Sons.\\

\nin Hardy, G.\ H.\ and Ramanujan, S.\ (1918).  Asymptotic formulae in combinatory analysis.  \ti{Proceedings of the London Mathematical Society} 17, 75--115.\\

\nin Hardy, G.\ H.\ and Wright, E.\ M.\ (1960).  \ti{The Theory of Numbers}, 4th ed.  Oxford:  Clarendon Press.\\

\nin  Hart, William L.\ (1948). Dunham Jackson 1888-1946. \ti{Bull. Amer. Math. Soc.} 54, 847--860.\\

\nin  Hart, William L.\ (1959). Dunham Jackson 1888-1946. \ti{Biographical Memoir}.  Washington, D.\ C.\.:  National Academy of Sciences.\\

\nin Hartman, Edwin P.\ (1970).  \ti{Adventures in Research:  A History of Ames Research Center 1940--1965}.  Washington, D.\ C.:  National Aeronautics and Space Administration.\\

\nin Heyde, C.\ C.\ and Seneta, E. (1977).  \ti{I.\ J.\ Bienaym\'{e}:  Statistical Theory Anticipated}.  New York:  Springer.\\

\nin Jackson, Bradley W.\ (n.\ d.).  A History of the Mathematics Department at San Jose State.  
\url{http://www.math.sjsu.edu/~jackson/} (accessed January 30, 2022).\\

\nin Jackson, Dunham (1935).  Mathematical principles in the theory of small samples.  \ti{The American Mathematical Monthly} 72, 344--364.\\  

\nin Jellison, Charles (1997).  The Prisoner and the Professor.  \ti{Stanford Magazine} March/April 1997.\\

\nin Kaplansky, Irving (2013).  Integers uniquely represented by certain ternary forms.  In R.L. Graham et al.\ (eds.), \ti{The Mathematics of Paul Erd\H{o}s I}.  New York:  Springer.\\

\nin Karatsuba, A.\ A.\ (1981).   Ivan Matveevich Vinogradov (on his ninetieth birthday).  \ti{Russian Mathematical Surveys}  36, 1--17.\\

\nin Krandick, Werner and Mehlhorn, Kurt (2006).  New bounds for the Descartes method.  \ti{Journal of Symbolic Computation} 41, 49--66.\\

\nin Krantz, Steven  (2019).  \ti{Mathematical Apocrypha Redux: More Stories and Anecdotes of Mathematicians and the Mathematical}.  American Mathematical Society, 83--84.\\

\nin Kuzmin, R.\ O.\  (1930). On a new class of transcendental numbers. \ti{Izvestiya Akademii Nauk SSSR}, Ser. matem. 7, 585--597.\\

\nin Lorentz, G.\ G.\ (2002).  Mathematics and politics in the Soviet Union from 1928 to 1953.  \ti{Journal of Approximation Theory} 116, 169--223.\\ 

\nin Malyshev, A. V.  and Faddeev, D. K. (1961).  Boris Alekseevich Venkov (on his sixtieth birthday). \ti{Uspekhi Mat. Nauk} 16, 235--240.\\

\nin Markov, A.\ A.\ (1888).  \ti{Table des valeurs de l'int\'{e}grale $\int_x^{\infty} e^{-t^2} dt$}. St.\ Petersburg:  Imperial Academy of Sciences.\\

\nin Markov, A.\ A.\ (1900).  \ti{Ischislenie Veroiatnostei} (Calculus of Probabilities).  Sanktpeterburg: Tipografija Imperatorskoj Akademii Nauk. Three subsequent editions:  2nd ed.\ 1908 and 3rd ed.\ 1913, Sanktpeterburg;  4th (posthumous) ed.\ 1924, Moskva: Gosudarstvennoe Izdatel'stvo. German translation of 3rd ed.:  \ti{Wahrscheinlichkeitsrechnung} (Heinrich Liebmann, trans.),  Leipzig and Berlin:  Teubner, 1912.\\

\nin McKay, Brendan (1989).  On Littlewood's estimate for the binomial distribution.  \ti{Advances in Applied Probability} 21, 475--478.\\

\nin Miller, G.\ H.\ (1970). "Blichfeldt, Hans Frederick", in \ti{Dictionary of Scientific Biography}, Volume 2.  New York:  Charles Scribner's Sons.  Reprinted 2008 in the ebook \ti{Complete Dictionary of Scientific Biography}.\\

\nin Munroe, M.\ E.\ (1951).  \ti{Theory of Probability}.  New York:  McGraw-Hill.\\

\nin Nazarov, Alexander I.\ and Sinkevich, Galina I.\ (2019).  History of Leningrad mathematics in the first half of the 20th century.  \ti{Notices of the international Congress of Chinese mathematicians}  7 (Number 2), 47--63.\\

\nin Nevai, Paul (1986).   G\'{e}za Freud, orthogonal polynomials, and Christoffel functions.  A case study.  \ti{Journal of Approximation Theory} 48, 3--167.\\

\nin Olds, C. D. (1963).  \ti{Continued Fractions}.  New Mathematical Library 9, Mathematical Association of America.  NY:  Random House.\\

\nin Olds, C. D. (1970).  The simple continued fraction expansion of $e$.  \ti{The American Mathematical Monthly} 77, 968--974.\\

\nin Ondar, Kh.\ O., ed. (1981).  \ti{The Correspondence Between A. \ A.\ Markov and A.\ A.\ Chuprov on the Theory of Probability and Mathematical Statistics} (translated by Charles and Margaret Stein).  New York:  Springer-Verlag.\\ 

\nin Ostrowski, A.\ M.\ (1950).  Note on Vincent's theorem.  \ti{Annals of Mathematics} 52, 702--707.\\

\nin Pearson, Egon (1990).  \ti{`Student':  A Statistical Biography of William Sealy Gosset} (R.\ L.\ Plackett and G.\ A.. Barnard, eds.).  Oxford:  Clarendon Press.\\

\nin Pearson, Karl (1931).  Historical note on the distribution of the standard deviations of samples of any size drawn from an indefinitely large normal parent population.  \ti{Biometrika} 23, 416--418.\\

\nin Peizer, David B.\ and Pratt, John W.\ (1968).  A normal approximation for binomial, F, beta, and other common, related tail probabilities, I.  \ti{Journal of the American Statistical Association} 63, 1416--1456.\\

\nin Perlman, Michael D.\ and Wichura, Michael J.\ (1975).  Sharpening Buffon's needle.  \ti{The American Statistician} 29, 157--163.\\

\nin P\'{o}lya, G., Szeg\"{o}, G.\, and Young, D.\ H.\ (1947).  Memorial Resolution:  James Victor Uspensky.  Stanford University Archives.\\

\nin Rashdall, Hastings (1895).  \ti{The Universities of Europe in the Middle Ages.  Volume 1:  Salerno---Bologna---Paris}.  Oxford:  Clarendon Press.\\

\nin Reeds, J.\ A.\, Diffie, W.\ et al.\ (2015).  \ti{Breaking Teleprinter Ciphers at Bletchley Park}.  	Piscataway, NJ:  Wiley--IEEE Press.\\

\nin Rietz, H.\ L.\ (ed.) (1924).  \ti{Handbook of Mathematical Statistics}.  Boston:  Houghton Mifflin.  \\

\nin Rouch\'{e}, E.\ (1888).  Sur un probl\`{e}me relatif \`{a} la dur\'{e}e du jeu.  \ti{C.\ R.\ Acad.\ Sci.} 116, 47--49.\\

\nin Royden, Halsey (1989).  The history of the Mathematics Department at Stanford.  In \ti{A Century of Mathematics in America, Part II} (ed.\ Peter Duren).  Providence, RI:  American Mathematical Society,  pp.\ 237--278.\\

\nin Royden, Halsey et al.\ (1992).  Memorial resolution:  Harold M. Bacon (1907--1992).  Stanford University Libraries, Department of Special Collections and University Archives.  Available at \url{https://purl.stanford.edu/ws360gk0273} (accessed February 29, 2020).\\

\nin Sanders, J.\ Thomas (1993).  The third opponent: dissertation defenses and the public profile of academic history in late Imperial Russia.  \ti{Jahrb\"{u}cher f\"{u}r Geschichte Osteuropas}, Neue Folge, 41, 242--265.\\

\nin Sansone, Giovanni (1950).  La formula di approssimazione asintotica dei polinomi di Tchebychef-Laguerre col procedimento di J.\ V.\ Uspensky.  \ti{Mathematische Zeitschrift} 53, 97--105.\\

\nin Seal, Hillary (1966).  The random walk of a simple risk business.  \ti{ASTIN Bulletin:  The Journal of the IAA} 4, 19--28.\\

\nin Siegmund-Schultze, Reinhard (2006).  Probability in 1919/20:  the von Mises-P\'{o}lya controversy.  \ti{Archive for History of Exact Sciences} 60, 431--515.\\

\nin Seneta, Eugene  (2006).  Markov and the the creation of Markov chains. In Amy N.
Langville and William J. Stewart (eds.), \ti{MAM2006: Markov Anniversary Meeting}.
Boson Books, Raleigh, North Carolina, pp. 1--20.\\

\nin Seneta, E. (2013).  A Tricentenary history of the Law of Large Numbers.  \ti{Bernoulli} 19,1088--1121.\\

\nin Sheynin, O.\ (2005).  \ti{Jakob Bernoulli:  On the Law of Large Numbers}, translated into English by Oscar Sheynin.  \url{http://sheynin.de/download/bernoulli.pdf} (last accessed December 20, 2021).\\

\nin Stanford (1992).  Long-time calculus professor Harold Bacon dies at 85.  Stanford University News Service.\\

\nin Stigler, Stephen M.\ (1986).  \ti{The History of Statistics:  The Measurement of Uncertainty before 1900}.  Cambridge, MA:  Harvard University Press.\\

\nin Strutt, John William, Baron Rayleigh (1894).  \ti{The Theory of Sound}, Vol.\ 1 (2nd ed.).  London and New York:  Macmillan.\\

\nin Sviderskaya, Galina Evgenievna (nee Kuzmina) (2002).  Rodion Osievich Kuzmin.  \ti{Public Scientific and Methodological Internet Journal of St. Petersburg State Technical University} - 2001--2002. - No. 2.\\

\nin Timoshenko, S.\ (1968)  \ti{As I Remember: The Autobiography of Stephen P. Timoshenko}, translated from the Russian by Robert Addis. Princeton, N.J.: Van Nostrand.\\

\nin Timoshenko, S.\ (1947).  Unpublished typescript (presumably prepared shortly after Uspensky's death and in the Stanford University Archives).\\

\nin Todhunter, I.\ (1865).  \ti{A History of the Mathematical Theory of Probability}.  London:  Macmillan.  Reprinted, 1949 and 1964, New York:  Chelsea.  [``A work of true learning, beyond criticism''---John Maynard Keynes.]\\

\nin Venkov, B.\ A.\, and Natanson, I.\ P.\, (1949).  Rodion Osievich Kuzmin (1891--1949) (obituary). \ti{Uspekhi Matematicheskikh Nauk} 4, 148--155.\\

\nin Verger, Jacques (2003).  Teachers.  In Ridder-Symoens, Hilde de (ed.), \ti{A History of the University in Europe.  Volume 1:  Universities in the Middle Ages}, Chapter 5.  Cambridge:  Cambridge University Press, pp.\ 144--168.\\

\nin Vinogradov, I.\ M.\ (1936).  \ti{Foundations of the Theory of Numbers}.  Moscow-Leningrad.  Eight editions, 1936--1971.  English translation of the fifth (1949) ed., \ti{Elements of Number Theory}, 1954 (trans.\ S.\ Kravetz), New York:  Dover Publications.\\

\nin Wall, H.\ (1948).  \ti{Analytic Theory of Continued Fractions}.  New York:  Chelsea.\\

\nin Wythoff, W.\ A.\ (1907).  A modification of the game of Nim.  \ti{Nieuw Archief voor Wiskunde} 7, 199-202.\\

\nin Zabell, S.\ L.\ (2008).  On Student's 1908 article ``The probable error of a mean''.  \ti{Journal of the American Statistical Association} 103, 1--7.\\

\end{document}